%% file: main.tex
\documentclass[letterpaper, 10 pt, conference]{ieeeconf}
\IEEEoverridecommandlockouts
\let\labelindent\relax
\usepackage{enumitem}
\usepackage{balance}
\usepackage[font={small}]{caption}
\usepackage{subcaption}
\usepackage{array}
\usepackage{textcomp}
\usepackage{mathtools, nccmath}
\usepackage{graphicx}
\usepackage{amsfonts}
\usepackage{amsmath}
\usepackage{amssymb}
\usepackage{algorithm}
\usepackage{algorithmic}
\usepackage{hyperref}
\usepackage{tikz}
\usepackage{arydshln}
\usepackage{multirow}
\usepackage{bm}
\usepackage{epstopdf}
\usepackage{cite}
\usepackage{siunitx}

\usepackage{amsthm}

\newtheorem{theorem}{Theorem}

\newtheorem{remark}{Remark}
\theoremstyle{definition}
\newtheorem{definition}{Definition}

\title{\LARGE \bf
Iterative Convex Optimization for Model Predictive Control \\ with Discrete-Time High-Order Control Barrier Functions
}

\author{Shuo Liu$^{*1}$, Jun Zeng$^{*2}$, Koushil Sreenath$^{2}$ and Calin A. Belta$^{1}$
\thanks{$^*$ Authors contributed equally.}
\thanks{This work was supported in part by the NSF under grants IIS-2024606 and CMMI-1931853.}
\thanks{$^{1}$S. Liu and C. Belta are with the department of Mechanical Engineering, Boston
University, Brookline, MA, 02215, USA {\tt\small \{liushuo, cbelta\}@bu.edu}. $^{2}$J. Zeng and K. Sreenath are with the University of California, Berkeley, CA, 94720, USA
{\tt\small  \{zengjunsjtu, koushils\}@berkeley.edu}}
\thanks{Implementation code is released on \url{https://github.com/ShockLeo/Iterative-MPC-DHOCBF}.}%
}

\begin{document} 
\maketitle
\input{sections/abstract.tex}
\input{sections/introduction.tex}
\input{sections/background.tex}

\input{sections/formulation.tex}
\input{sections/results.tex}
\input{sections/conclusion.tex}

\bibliographystyle{IEEEtran}
\balance
\bibliography{references.bib}

\end{document}

%% file: sections/abstract.tex
\begin{abstract}
Safety is one of the fundamental challenges in control theory. 
Recently, multi-step optimal control problems for discrete-time dynamical systems were formulated to enforce stability, while subject to input constraints as well as safety-critical requirements using discrete-time control barrier functions within a model predictive control (MPC) framework.
Existing work usually focus on the feasibility or the safety for the optimization problem, and the majority of the existing work restrict the discussions to relative-degree one control barrier functions.
Additionally, the real-time computation is challenging when a large horizon is considered in the MPC problem for relative-degree one or high-order control barrier functions.
In this paper, we propose a framework that solves the safety-critical MPC problem in an iterative optimization, which is applicable for any relative-degree control barrier functions.
In the proposed formulation, the nonlinear system dynamics as well as the safety constraints modeled as discrete-time high-order control barrier functions (DHOCBF) are linearized at each time step.
Our formulation is generally valid for any control barrier function with an arbitrary relative-degree.
The advantages of fast computational performance with safety guarantee are analyzed and validated with numerical results. 
\end{abstract}

%% file: sections/introduction.tex
\section{Introduction}
\subsection{Motivation}
Safety-critical optimal control is a central problem in robotics.
For example, reaching a goal while avoiding obstacles and minimizing energy can be formulated as a constrained optimal control problem by using continuous-time control barrier functions (CBFs)~\cite{ames2014control, ames2016control}. By dividing the timeline into small intervals, the problem is reduced to a (possibly large) number of  quadratic programs, which can be solved at real-time speeds. However, this approach can be too aggressive due to the lack of predicting ahead.

Model predictive control (MPC) with CBFs~\cite{zeng2021safety} considers the safety problem in the discrete-time domain, and provides a smooth control policy as it involves future state information along a receding horizon.
However, the computational time is relatively large and increases dramatically with a larger horizon, since the optimization itself is usually nonlinear and non-convex.
An additional issue of this nonlinear model predictive formulation is the feasibility of the optimization. For CBFs with relative-degree one, relaxation techniques have been introduced in~\cite{zeng2021enhancing}.
In this paper, we address the above challenges with a proposed convex MPC with linearized, discrete-time CBFs, under an iterative approach. In contrast with the real-time iteration (RTI) approach introduced in ~\cite{diehl2005real}, which solves the problem through iterative Newton steps, our approach solves the optimization problem formulated by a convex MPC iteratively for each time step.
We show that the proposed approach can significantly reduce the computational time, compared to the state of the art introduced in~\cite{zeng2021enhancing}, even for CBFs with high relative-degree, without sacrificing the controller performance.
The feasibllity rate of our proposed method also outperforms that of the baseline method in~\cite{zeng2021enhancing} for large horizon lengths.

\subsection{Related work}

\subsubsection{Model Predictive Control (MPC)}
MPC is widely used in modern control systems, such as controller design in robotic manipulation and locomotion~\cite{paolo2017mpc, scianca2020mpc} to obtain a control strategy as a solution to an optimization problem. Stability was achieved in \cite{grandia2020nonlinear} by incorporating discrete-time control Lyapunov functions (DCLFs) into a general MPC-based optimization problem to realize real-time control on a robotic system with limited computational resources. More and more recent work like \cite{son2019safety} emphasizes safety in robot design and deployment since it is an important criterion for real-world tasks. Some works consider safety criteria through the introduction of additional repelling functions~\cite{ames2014control, eklund2011switched} while some works  regard obstacle avoidance as one concrete scenario in terms of safety criteria for robots \cite{frasch2013auto, liniger2015optimization, zhang2020optimization}. Those safety criteria are usually formulated as constraints in optimization problems.
 This paper can be seen in the context of MPC with safety constraints.

\subsubsection{Continuous-Time CBFs}
It has recently been shown that to stabilize an affine control system while also satisfying safety constraints and control limitations, CBFs can be unified with control Lyapunov functions (CLFs) to form a sequence of single-step optimization programs~\cite{ames2014control, ames2016control, galloway2015torque, ames2019control}.
If the cost is quadratic, the optimizations are quadratic programs (QP), and the solutions can be deployed in real time \cite{ames2014control, nguyen20163d}.
Adaptive, robust and stochastic versions of safety-critical control with CBFs were introduced in~\cite{xiao2021adaptive, Nguyen2015_RobustCLF, jankovic2018robust, clark2021control, dawson2022safe}. 
For safety constraints expressed using functions with high relative degree with respect to the dynamics of the system, exponential CBFs~\cite{nguyen2016exponential} and high-order CBFs (HOCBFs)~\cite{xiao2021high, tan2021high, usevitch2022adversarial} were proposed.

\subsubsection{Discrete-Time CBFs}
Discrete-time CBFs (DCBFs) were introduced in \cite{agrawal2017discrete} as a means to enable safety-critical control for discrete-time systems. They were used in a nonlinear MPC (NMPC) framework to create NMPC-DCBF \cite{zeng2021safety}, wherein the DCBF constraint was enforced through a predictive horizon.
This method was also utilised in a multi-layer control framework in~\cite{grandia2021multi}, where DCBFs with longer horizons were considered in the MPC problem serving as a mid-level controller to guarantee safety.
Generalized discrete-time CBFs (GCBFs) and discrete-time high-order CBFs (DHOCBFs) were proposed in \cite{ma2021feasibility} and \cite{xiong2022discrete} respectively, where the DCBF constraint only acted on the first time-step, i.e., a single-step constraint.
MPC with DCBF has been used in various fields, such as autonomous driving~\cite{he2022autonomous} and legged robotics~\cite{li2022bridging}.
For the work above, the CBF constraints are either limited to be activated at the first time-step~\cite{agrawal2017discrete, ma2021feasibility, xiong2022discrete} to improve the optimization feasibility at the cost of sacrificing the safety performance, or for multiple or all steps~\cite{grandia2021multi, he2022autonomous, li2022bridging} with additional performance optimization from other modules, such as multi-layer control~\cite{grandia2021multi, he2022autonomous} or planning~\cite{li2022bridging}, which needs to specified for different platforms.
A decay-rate relaxing technique~\cite{zeng2021decay} was introduced for NMPC with DCBF~\cite{zeng2021enhancing} for all time-steps to enhance the safety and feasibility at the same time, but the computation itself is overall still nonlinear and non-convex which could be computationally slow for large horizons and nonlinear dynamical systems, and the discussion in \cite{zeng2021enhancing} is limited to relative-degree one.
In this paper, we generalize relaxing technique for DHOCBF and largely optimize the computational time compared to all existing work.

\subsection{Contributions}

We propose a novel approach to the NMPC with discrete-time CBFs that is significantly faster than existing approaches. In particular, the contributions are as follows:

\begin{itemize}
\item We present a model predictive control strategy for safety-critical tasks, where the safety-critical constraints can be enforced by DHOCBFs. The decay rate in each constraint can be relaxed to enhance the feasibility in optimization and to ensure forward invariance of the intersection of a series of safety sets.
\item We propose an optimal control framework for guaranteeing safety, where the DHOCBF constraints as well as the system dynamics are linearized at each iteration, and considered as constraints in a convex optimization solved iteratively.
\item We show through numerical examples that the proposed framework is significantly faster than existing methods, without sacrificing safety and feasibility.
\end{itemize}

%% file: sections/background.tex
  \section{Preliminaries}

In this section, we introduce some definitions and results on CBF and MPC.

\subsection{Discrete-Time High-Order Control Barrier Function (DHOCBF)}
\label{subsec:dhocbf}
In this work, safety is defined as forward invariance of a set $\mathcal C$, i.e., a system is said to be {\em safe} if it stays in $\mathcal C$ for all time, given that it is initialized in $\mathcal C$. We consider the set  $\mathcal C$ as the superlevel set of a discrete-time function $h: \mathbb{R}^{n}\to \mathbb{R}$:
\begin{equation}
\label{eq:safe-set}
\mathcal{C} \coloneqq \{\mathbf{x}_{t} \in \mathbb{R}^{n}: h(\mathbf{x}_{t} )\geq0 \}.
\end{equation}
We consider a discrete-time control system in the form  
\begin{equation}
\label{eq:discrete-dynamics}
\mathbf{x}_{t+1}=f(\mathbf{x}_{t},\mathbf{u}_{t}),
\end{equation}
where $\mathbf{x}_{t} \in \mathcal X \subset \mathbb{R}^{n}$ represents the state of system \eqref{eq:discrete-dynamics} at time step $t\in\mathbb{N}, \mathbf{u}_{t}\in \mathcal U \subset \mathbb{R}^{q}$ is the control input, and function $f$ is locally Lipschitz.

\begin{definition}[Relative degree~\cite{sun2003initial}]
\label{def:relative-degree}
The output $\mathbf{y}_{t}=h(\mathbf{x}_{t})$ of system \eqref{eq:discrete-dynamics} is said to have relative degree $m$ if
\begin{equation}
\begin{split}
\mathbf{y}_{t+i}=h(\bar{f}_{i-1}(f(\mathbf{x}_{t},\mathbf{u}_{t}))), \ i \in \{1,2,\dots,m\},\\
\text{s.t.} \ \frac{\partial \mathbf{y}_{t+m}}{\partial \mathbf{u}_{t}}\ne 0_{q}, \frac{\partial \mathbf{y}_{t+i}}{\partial \mathbf{u}_{t}}=0_{q},  \ i \in \{1,2,\dots,m-1\},
\end{split}
\end{equation}
i.e., $m$ is the number of steps (delay) in the output $\mathbf{y}_{t}$ in order for the control input $\mathbf{u}_{t}$ to appear. 
\end{definition}

In the above definition, we use $\bar{f}(\mathbf{x}_{t})$ to denote the uncontrolled state dynamics $f(\mathbf{x}_{t}, 0)$. The subscript $i$ of function $\bar{f}(\cdot)$ denotes the $i$-times recursive compositions of $\bar{f}(\cdot)$, i.e.,  $\bar{f}_{i}(\mathbf{x}_{t})=\underset{i\text{-times}~~~~~~~~~~~~~}{\underbrace{\bar{f}(\bar{f}(\dots,\bar{f}}(\bar{f}_{0}(\mathbf{x}_{t}))))}$ with $\bar{f}_{0}(\mathbf{x}_{t})=\mathbf{x}_{t}$.

We assume that $h(\mathbf{x}_{t})$ has
relative degree $m$ with respect to system (\ref{eq:discrete-dynamics}) based on Def. \ref{def:relative-degree}.
Starting with $\psi_{0}(\mathbf{x}_{t})\coloneqq h(\mathbf{x}_{t})$, we define a sequence of discrete-time functions $\psi_{i}:  \mathbb{R}^{n}\to\mathbb{R}$, $i=1,\dots,m$ as:
\begin{equation}
\label{eq:high-order-discrete-CBFs}
\psi_{i}(\mathbf{x}_{t})\coloneqq \bigtriangleup \psi_{i-1}(\mathbf{x}_{t},\mathbf{u}_{t})+\alpha_{i}(\psi_{i-1}(\mathbf{x}_{t})), 
\end{equation}
where $\bigtriangleup \psi_{i-1}(\mathbf{x}_{t}, \mathbf{u}_{t})\coloneqq \psi_{i-1}(\mathbf{x}_{t+1})-\psi_{i-1}(\mathbf{x}_{t})$, and $\alpha_{i}(\cdot)$ denotes the $i^{th}$ class $\kappa$ function which satisfies $\alpha_{i}(\psi_{i-1}(\mathbf{x}_{t}))\le \psi_{i-1}(\mathbf{x}_{t})$ for $i=1,\ldots, m$.
A sequence of sets $\mathcal {C}_{i}$ is defined based on \eqref{eq:high-order-discrete-CBFs} as
\begin{equation}
\label{eq:high-order-safety-sets}
\mathcal {C}_{i}\coloneqq \{\mathbf{x}_{t}\in \mathbb{R}^{n}:\psi_{i}(\mathbf{x}_{t})\ge 0\}, \ i =\{0,\ldots,m-1\}.
\end{equation}

\begin{definition}[DHOCBF~\cite{xiong2022discrete}]
\label{def:high-order-discrete-CBFs}
Let $\psi_{i}(\mathbf{x}_{t}), \ i\in \{1,\dots,m\}$ be defined by \eqref{eq:high-order-discrete-CBFs} and $\mathcal {C}_{i},\ i\in \{0,\dots,m-1\}$ be defined by \eqref{eq:high-order-safety-sets}. A function $h:\mathbb{R}^{n}\to\mathbb{R}$ is a Discrete-Time High-Order Control Barrier Function (DHOCBF) with relative degree $m$ for system \eqref{eq:discrete-dynamics} if there exist $\psi_{m}(\mathbf{x}_{t})$ and $\mathcal {C}_{i}$ such that
\begin{equation}
\label{eq:highest-order-CBF}
\psi_{m}(\mathbf{x}_{t})\ge 0, \ \forall x_{t}\in \mathcal{C}_{0}\cap \dots \cap \mathcal {C}_{m-1}.
\end{equation}
\end{definition}

\begin{theorem}[Safety Guarantee \cite{xiong2022discrete}]
\label{thm:forward-invariance}
Given a DHOCBF $h(\mathbf{x}_{t})$ from Def. \ref{def:high-order-discrete-CBFs} with corresponding sets $\mathcal{C}_{0}, \dots,\mathcal {C}_{m-1}$ defined by \eqref{eq:high-order-safety-sets}, if $\mathbf{x}_{0} \in \mathcal {C}_{0}\cap \dots \cap \mathcal {C}_{m-1},$ then any Lipschitz controller $\mathbf{u}_{t}$ that satisfies the constraint in \eqref{eq:highest-order-CBF}, $\forall t\ge 0$ renders $\mathcal {C}_{0}\cap \dots \cap \mathcal {C}_{m-1}$ forward invariant for system \eqref{eq:discrete-dynamics}, $i.e., \mathbf{x}_{t} \in \mathcal {C}_{0}\cap \dots \cap \mathcal {C}_{m-1}, \forall t\ge 0.$
\end{theorem}

\begin{remark}
\label{rem:sufficient-condition}
The function $\psi_{i}(\mathbf{x}_{t})$ in \eqref{eq:high-order-discrete-CBFs} is called a $i^{th}$ order discrete-time control barrier function (DCBF) in this paper. Since satisfying the $i^{th}$ order DCBF constraint ($\psi_{i}(\mathbf{x}_{t})\ge0$) is a sufficient condition for rendering $\mathcal{C}_{0}\cap \dots \cap \mathcal{C}_{i-1}$ forward invariant for system \eqref{eq:discrete-dynamics} as shown above, it is not necessary to formulate DCBF constraints up to $m^{th}$ order as \eqref{eq:highest-order-CBF} if the control input $\mathbf{u}_{t}$ could be involved in some optimal control problem, which allows us to choose an appropriate order for the DCBF constraint to reduce the computation. In other words, the highest order for DCBF could be $m_{\text{cbf}}$ with $m_{\text{cbf}}\le m$. We can simply define a $i^{th}$ order DCBF $\psi_{i}(\mathbf{x}_{t})$ in \eqref{eq:high-order-discrete-CBFs} as
\begin{equation}
\label{eq:simple-high-order-discrete-CBFs}
\psi_{i}(\mathbf{x}_{t})\coloneqq \bigtriangleup \psi_{i-1}(\mathbf{x}_{t},\mathbf{u}_{t})+\gamma_{i}\psi_{i-1}(\mathbf{x}_{t}),
\end{equation}
where $0<\gamma_{i}\le 1, i\in \{1,\dots,m_{\text{cbf}}\}$.
\end{remark}

The expression in~\eqref{eq:simple-high-order-discrete-CBFs} follows the format of the first order DCBF proposed in \cite{agrawal2017discrete} and could be used to define a DHOCBF with arbitrary relative degree. 

\subsection{Model Predictive Control}

Consider the problem of regulating to a target state $\mathbf{x}_{r}$ for the discrete-time system \eqref{eq:discrete-dynamics} while making sure that safety is guaranteed by ensuring $\psi_{0}(\mathbf{x}_{t})=h(\mathbf{x}_{t})\ge 0$. 
The following optimal control problem takes future $N$ states into account as prediction at each time step $t$:

\noindent\rule{\columnwidth}{0.4pt}
\textbf{NMPC-DCBF:}
{\small
\begin{subequations}
\label{eq:mpc-cbf}
\begin{align}
\label{subeq:mpc-cbf-cost}
\min_{\mathbf{U}_{t},\Omega_{t}} \ & p(\mathbf{x}_{t,N})+\sum_{k=0}^{N-1}q(\mathbf{x}_{t,k},\mathbf{u}_{t,k},\omega_{t,k}) \\
\text{s.t.} \ & \mathbf{x}_{t, k+1} = f(\mathbf{x}_{t, k}, \mathbf{u}_{t, k}), \ k{=}\{0,\dots,N{-}1\} \label{subeq:mpc-cbf-dyn} \\
& \mathbf{u}_{t, k} \in \mathcal U, \mathbf{x}_{t, k} \in \mathcal X, \omega_{t, k}\in \mathbb{R},\ k{=}\{0,\dots,N{-}1\} \label{subeq:mpc-cbf-cons1} \\
& h(\mathbf{x}_{t, k+1}) \ge \omega_{t, k}(1-\gamma)h(\mathbf{x}_{t, k}), 0<\gamma\le1, \label{subeq:mpc-cbf-cons2} \\ 
& k=\{0,\dots, N{-}1\}, \nonumber
\end{align}
\end{subequations}
}
\noindent\rule{\columnwidth}{0.4pt}
where $\mathbf{x}_{t, k+1}$ denotes the state at time step $k+1$ predicted at time step $t$ obtained by applying the input vector $\mathbf{u}_{t, k}$ to the state $\mathbf{x}_{t, k}$.
In \eqref{subeq:mpc-cbf-cost}, $q(\cdot)$ and $p(\cdot)$ denote stage and terminal costs, respectively, and $\omega_{t,k}$ is a slack variable.
The discrete-time dynamics is represented by \eqref{subeq:mpc-cbf-dyn} and the constraints of state and control input along the horizon are captured by \eqref{subeq:mpc-cbf-cons1}.
The DCBF constraint in \eqref{subeq:mpc-cbf-cons2} is proposed in \cite{agrawal2017discrete} and is designed to ensure the forward invariance of the set $\mathcal{C}$ based on \eqref{eq:safe-set}.
The above formulation was first proposed in~\cite{zeng2021safety}, and then later generalized in~\cite{zeng2021enhancing}, where the decay rate $(1-\gamma)$ of the CBF was relaxed by slack variable $\omega_{t,k}$ to enhance safety and feasibility.

The optimal solution to \eqref{eq:mpc-cbf} at time $t$ is $\mathbf{U}^{*}_{t}=[\mathbf{u}_{t, 0}^{*},\dots,\mathbf{u}_{t,N-1}^{*}]$ and $\Omega^{*}_{t}=[\omega_{t, 0}^{*},\dots,\omega_{t, N-1}^{*}]$. The first element of $\mathbf{U}_{t}^{*}$ is applied to \eqref{eq:discrete-dynamics} as
\begin{equation}
\label{eq:opt-trajectory}
\mathbf{x}_{t+1}=f(\mathbf{x}_{t},\mathbf{u}_{t,0}^{*})
\end{equation}
to get the new state $\mathbf{x}_{t+1}$. The constrained finite-time optimal control problem \eqref{eq:mpc-cbf} is solved at time step $t+1$, and all future time steps based on the new state $\mathbf{x}_{t+1}$, yielding a safety-critical receding horizon control strategy.

%% file: sections/formulation.tex
\section{Iterative Convex MPC with DHOCBF}
\label{sec:iterative-opt}
In this section, we present an iterative convex MPC for DCBF, which works for general DHOCBFs defined in Sec. \ref{subsec:dhocbf}.

\begin{algorithm}
\caption{iMPC-DHOCBF}
\begin{algorithmic}[1]
\label{alg:iMPC-DCBF}
\begin{small}
 \renewcommand{\algorithmicrequire}{\textbf{Input:}}
 \renewcommand{\algorithmicensure}{\textbf{Output:}}
 \REQUIRE System dynamics \eqref{eq:discrete-dynamics}, candidate CBF constraint, obstacle configurations, initial state $\mathbf{x}(0)$.
 \ENSURE Safety-critical optimal control for obstacle avoidance. \\
 \STATE 
 Set initial guess $\bar{\mathbf{U}}_{0}^{0} = 0$ at $t = 0$.
 \STATE
 Propagate with system dynamics to get initial guess of states $\bar{\mathbf{X}}_{0}^{0}$ from initial state $\mathbf{x}(0)$.
    \FOR{$t\le t_{\text{sim}}-1$}
        \STATE
        Initialize $j = 0$.
        \WHILE{ Iteration $j$ (not converged OR $j < j_{\text{max}}$)} 
          \STATE Linearize system dynamics / constraints with $\bar{\mathbf{X}}_{t}^{j}, \bar{\mathbf{U}}_{t}^{j}$.
          \STATE
          Solve a convex finite-time constrained optimal control problem (CFTOC) with linearized dynamics / constraints and get optimal values of states and inputs $\mathbf{X}_{t}^{*,j}$, $\mathbf{U}_{t}^{*,j}$.
          \STATE
          $\bar{\mathbf{X}}_{t}^{j+1} = \mathbf{X}_{t}^{*, j}, \bar{\mathbf{U}}_{t}^{j+1} = \mathbf{U}_{t}^{*, j}$,
          $j = j+1$
        \ENDWHILE
        \STATE
        Extract optimized states and inputs $\mathbf{X}_{t}^{*} = \mathbf{X}_{t}^{*,j}, \mathbf{U}_{t}^{*} = \mathbf{U}_{t}^{*,j}$ from last iteration and extract $\mathbf{u}_{t,0}^{*}$ from $\mathbf{U}_{t}^{*}$.
        \STATE 
        Apply $\mathbf{u}_{t,0}^{*}$ with respect to system dynamics \eqref{eq:discrete-dynamics} to get $\mathbf{x}_{t+1} = f(\mathbf{x}_{t}, \mathbf{u}_{t,0}^{*})$, and record $\mathbf{x}(t+1) = \mathbf{x}_{t+1}$.
        \STATE
        Update $\bar{\mathbf{U}}_{t+1}^{0}$ with $\mathbf{U}_{t}^{*}$ and propagate to calculate $\bar{\mathbf{X}}_{t+1}^{0}$.
        \STATE
        $t = t+1$.
    \ENDFOR
 \RETURN closed-loop trajectory $[\mathbf{x}(0),\dots,\mathbf{x}(t_{\text{sim}})]$
\end{small}
\end{algorithmic}
\end{algorithm}

\subsection{Iterative Convex Optimization}
\label{subsec:iterative-convex-optimization}

The algorithm described in Alg.~\ref{alg:iMPC-DCBF}
contains an iterative optimization at each time step $t$, which is denoted as \emph{iterative} MPC-DHOCBF (iMPC-DHOCBF).  
Our iterative optimization problem contains three parts for each iteration $j$: (1) solve a convex finite-time optimal control (CFTOC) problem with linearized dynamics and DHOCBF, (2) check convergence criteria, (3) update state and input vectors for next iteration.
Notice that the open-loop trajectory with updated states $\bar{\mathbf{X}}_{t}^{j} = [\bar{\mathbf{x}}_{t, 0}^{j},\dots, \bar{\mathbf{x}}_{t, N-1}^{j}]$ and inputs $\bar{\mathbf{U}}_{t}^{j} = [\bar{\mathbf{u}}_{t, 0}^{j},\dots, \bar{\mathbf{u}}_{t, N-1}^{j}]$ is passed between iterations, which allows iterative linearization for both system dynamics and DHOCBF locally.
As discussed before, ``high-order" implies that the relative degree should be larger or equal to one.

The iteration is finished when the convergence error function $e(\mathbf{X}_{t}^{*, j}, \mathbf{U}_{t}^{*, j}, \bar{\mathbf{X}}_{t}^{j}, \bar{\mathbf{U}}_{t}^{j})$ is within a user-defined normalized convergence criteria, where $\mathbf{X}_{t}^{*, j}=[\mathbf{x}_{t, 0}^{*, j},\dots, \mathbf{x}_{t, N}^{*, j}]$, $\mathbf{U}_{t}^{*, j}=[\mathbf{u}_{t,0}^{*, j},\dots, \mathbf{u}_{t,N-1}^{*, j}]$ represent optimized states and inputs at iteration $j$.
To restrict the number of iterations, we limit $j < j_{\text{max}}$, where $j_{\text{max}}$ denotes the maximum numbers of iterations.
Therefore, the iterative optimization stops when the cost function reaches a local optimal minimum, whose iteration number is denoted as $j_{t, \text{conv}}$.
The optimized states $\mathbf{X}_{t}^{*}$ and inputs $\mathbf{U}_{t}^{*}$ are passed to the iMPC-DHOCBF formulation for the next time instant. 
At each time, we record the updated state propagated by the system dynamics with a given discretization time, which allows to extract the output closed-loop trajectory with our proposed iMPC-DHOCBF.
 
\subsection{Linearization of Dynamics}
\label{subsec:linearization-dynamics}
At iteration $j$, an improved vector $\mathbf{u}_{t,k}^{j}$ is considered by linearizing the system around $\bar{\mathbf{x}}_{t,k}^{j}, \bar{\mathbf{u}}_{t,k}^{j}$:
\begin{equation}
\begin{split}
\label{eq:linearized-dynamics}
\mathbf{x}_{t,k+1}^{j}{-}\bar{\mathbf{x}}_{t,k+1}^{j}{=}A^{j}(\mathbf{x}_{t,k}^{j}{-}\bar{\mathbf{x}}_{t,k}^{j})+B^{j}(\mathbf{u}_{t,k}^{j}{-}\bar{\mathbf{u}}_{t,k}^{j}),
\end{split}
\end{equation}
where $0 \leq j < j_{\text{max}}$; $k$ and $j$ represent open-loop time step and iteration indices, respectively. We also have
\begin{equation}
A^{j}=D_{\mathbf{x}}f(\bar{\mathbf{x}}_{t,k}^{j}, \bar{\mathbf{u}}_{t,k}^{j}), \ B^{j}=D_{\mathbf{u}}f(\bar{\mathbf{x}}_{t,k}^{j}, \bar{\mathbf{u}}_{t,k}^{j}),
\end{equation}
where $D_{\mathbf{x}}$ and $D_{\mathbf{u}}$ denote the Jacobian of the system dynamics $f(\mathbf{x}, \mathbf{u})$ with respect to the state $\mathbf{x}$ and the input $\mathbf{u}$.
This approach allows to linearize the system at $(\bar{\mathbf{x}}_{t,k}^{j}, \bar{\mathbf{u}}_{t,k}^{j})$ locally between iterations.
The convex system dynamics constraints are provided in ~\eqref{eq:linearized-dynamics} since all nominal vectors $(\bar{\mathbf{x}}_{t,k}^{j}, \bar{\mathbf{u}}_{t,k}^{j})$ in current iteration are constant and constructed from previous iteration $j-1$.

\subsection{Linearization of DCBF \& DHOCBF}
\label{subsec:linearization-dhocbf}

\begin{figure}
    \centering
    \includegraphics[width=\linewidth]{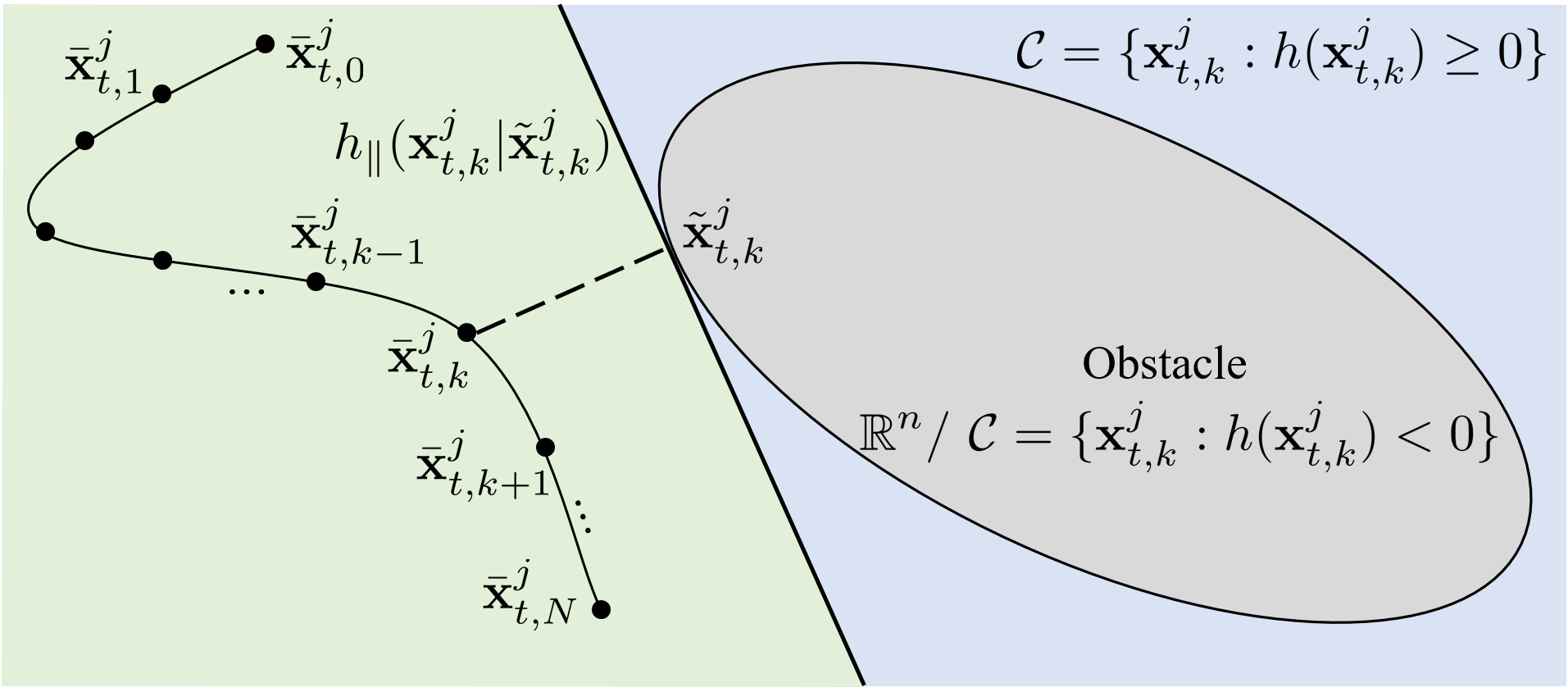}
    \caption{Linearization of DHOCBF: $h_{\parallel}(\mathbf{x}_{t, k}^{j}|\tilde{\mathbf{x}}_{t, k}^{j})\ge 0$ represents the linearized safe set locally and is colored in green. Note that $h_{\parallel}(\mathbf{x}_{t, k}^{j}|\tilde{\mathbf{x}}_{t, k}^{j})\ge 0$ guarantees $h(\mathbf{x}_{t,k}^{j}) \ge 0$ (colored in blue plus green), which ensures collision avoidance (outside the grey region).}
    \label{fig:linearization-dhocbf}
\end{figure}

In this section, we show how to linearize the DCBF up to the highest order. At iteration $j$, in order to linearize $h(\mathbf{x}_{t,k}^{j})$, an explicit line is projected in the state space to the nearest point $\tilde{\mathbf{x}}_{t, k}^{j}$ on the boundary of the obstacle from each state $\bar{\mathbf{x}}_{t, k}^{j}$. Note that $\bar{\mathbf{x}}_{t, k}^{j}$ is the nominal state vector from iteration $j-1$ for the linearization at iteration $j$, which means $\bar{\mathbf{x}}_{t, k}^{j}=\mathbf{x}_{t,k}^{j-1}$. The tangent line passing through the nearest point $\tilde{\mathbf{x}}_{t, k}^{j}$ is denoted as $h_{\parallel}(\mathbf{x}_{t, k}^{j}|\tilde{\mathbf{x}}_{t, k}^{j})$. 
This allows us to define a linearized safe set by $h_{\parallel}(\mathbf{x}_{t, k}^{j}|\tilde{\mathbf{x}}_{t, k}^{j})\ge 0$, $\forall t \in \mathbb{N}$ as shown in Fig. \ref{fig:linearization-dhocbf} by the green region.  

\begin{remark}
Note that $\tilde{\mathbf{x}}_{t, k}^{j}$ represents the optimized value of the minimum distance problem with distance function $h(\cdot)$ between $\bar{\mathbf{x}}_{t, k}^{j}$ and safe set $\mathcal{C}$.
For common smooth and differentiable CBFs, the expression of $\tilde{\mathbf{x}}_{t, k}^{j}$ as a function of $\bar{\mathbf{x}}_{t, k}^{j}$ is explicit~\cite{thirugnanam2022duality, thirugnanam2022safety}.
For example, when $h(\cdot)$ describes a $l_2$-norm function with the obstacle being a circular shape, $\tilde{\mathbf{x}}_{t, k}^{j}$ is exactly the intersection point between $\bar{\mathbf{x}}_{t, k}^{j}$ and the center of the obstacle.
Notice that $\tilde{\mathbf{x}}_{t, k}^{j}$ could be implicit for general elliptic calculations~\cite{skrypnik1994methods}, but it could still be numerically approximated as the values of $\bar{\mathbf{x}}_{t, k}^{j}$ known at iteration $j$ before the linearization.
\end{remark}

The relative degree of $h_{\parallel}(\mathbf{x}_{t,k}^{j}|\tilde{\mathbf{x}}_{t, k}^{j})$ with respect to system \eqref{eq:discrete-dynamics} is still $m$ when the relative degree of $h(\mathbf{x}_{t,k}^{j})$ is $m$.
Thus, in order to guarantee safety with forward invariance based on Thm.~\ref{thm:forward-invariance} and Rem.~\ref{rem:sufficient-condition}, two sufficient conditions need to be satisfied: (1) the sequence of linearized DHOCBF $\tilde{\psi}_{0}(\cdot),\dots, \tilde{\psi}_{m_{\text{cbf}}-1}(\cdot)$  is larger or equal to zero at the initial condition $\mathbf{x}_{t}$, and (2) the highest-order DCBF constraint $\tilde{\psi}_{m_{\text{cbf}}}(\mathbf{x}) \ge 0$ is always satisfied, where $\tilde{\psi}_{i}(\cdot)$ is defined as:
\begin{equation}
\label{eq:linearized-CBFs}
\begin{split}
 \tilde{\psi}_{0}(\mathbf{x}_{t,k}^{j}) \coloneqq & h_{\parallel}(\mathbf{x}_{t,k}^{j}|\tilde{\mathbf{x}}_{t, k}^{j}) \\
 \tilde{\psi}_{i}(\mathbf{x}_{t,k}^{j}) \coloneqq & \tilde{\psi}_{i-1}(\mathbf{x}_{t,k+1}^{j}){-}\tilde{\psi}_{i-1}(\mathbf{x}_{t,k}^{j}){+}\gamma_{i}\tilde{\psi}_{i-1}(\mathbf{x}_{t,k}^{j}). 
 \end{split}
\end{equation}
Here, we have $0 <\gamma_{i} \le 1$, $i\in\{1,\dots,m_{\text{cbf}}\}$, and $m_{\text{cbf}} \le m$ (as in \eqref{eq:simple-high-order-discrete-CBFs}).

\begin{remark}
From Rem.~\ref{rem:sufficient-condition}, it follows that that $m_{\text{cbf}}$ is not necessarily equal to $m$. A detailed discussion on this can be found in~\cite{ma2021feasibility, zeng2021enhancing}.
\end{remark}

An important issue is {\em feasibility}.
It is possible that $\psi_{i}(\mathbf{x}_{t,k}^{0}) \ge 0$, $1\le i \le m_{\text{cbf}}-1$, with $k\in\{0,\dots,N\}$ is not satisifed since the linearized DHOCBF functions $\tilde{\psi}_0(\cdot),\dots, \tilde{\psi}_{m_{\text{cbf}}-1}(\cdot)$ are more conservative than the original forms $\psi_0(\cdot),\dots, \psi_{m_{\text{cbf}}-1}(\cdot)$.
This problem can occur when the horizon is too large, or the linearization is too conservative. 
In order to handle this issue, we introduce a slack variable $\omega_{t,k,i}^{j}$ with a corresponding decay rate $(1-\gamma_{i})$:
\begin{equation}
\label{eq:nonconvex-hocbf-constraint}
  \tilde{\psi}_{i-1}(\mathbf{x}_{t,k+1}^{j})\ge \omega_{t,k,i}^{j}(1-\gamma_{i}) \tilde{\psi}_{i-1}(\mathbf{x}_{t,k}^{j}),\ \omega_{t,k,i}^{j} \in \mathbb{R},
\end{equation}
where $i\in\{1,\dots,m_{\text{cbf}}\}$.
The slack variable $\omega_{t,k,i}^{j}$ is selected by minimizing a cost function term to satisfy DCBF constraints at initial condition at any time step \cite{zeng2021enhancing}.

Another challenge induced by the DCBF linearization is that the constraints in \eqref{eq:nonconvex-hocbf-constraint} could be non-convex, since now $\omega_{t,k,i}^{j}$ and $\mathbf{x}_{t,k}^{j}$ are both optimization variables.
Note that $\tilde{\psi}_{0}(\mathbf{x}_{t,0}^{j})$ are constant, thus we can only place $\omega_{t,k,i}^{j}$ in front of $\tilde{\psi}_{0}(\mathbf{x}_{t,0}^{j})$ and move the other optimization variables to the other side of the inequalities.
This motivates us to provide the following form for reformulating \eqref{eq:nonconvex-hocbf-constraint} as convex constraints:
\begin{equation}
\label{eq:convex-hocbf-constraint}
\begin{split}
& \tilde{\psi}_{i-1}(\mathbf{x}_{t,k}^{j}) + \sum_{\nu =1}^{i}Z_{\nu,i}(1-\gamma_{i})^{k}\tilde{\psi}_{0}(\mathbf{x}_{t,\nu}^{j})\ge \\
& \quad \omega_{t,k,i}^{j} Z_{0,i}(1-\gamma_{i})^{k}\tilde{\psi}_{0}(\mathbf{x}_{t,0}^{j}), \\
& j \le j_{\text{max}} \in \mathbb{N}^{+}, \ i\in \{1,\dots,m_{\text{cbf}}\},  \ \omega_{t,k,i}^{j}\in  \mathbb{R}.
\end{split}
\end{equation}
In the above, $Z_{\nu,i}$ is a constant that can be obtained recursively by reformulating $\tilde{\psi}_{i-1}(\cdot)$ back to $\tilde{\psi}_{0}(\cdot)$ given $\nu \in\{0,..,i\}$.
We define $Z_{\nu,i}$ as follows. When $2\le i, \nu \le i-2$, we have
\begin{equation}
\label{eq:z-value}
\begin{split}
Z_{\nu,i}=\sum_{l=1}^{l_{\text{max}}}[(\gamma_{\zeta_{1}}-1)(\gamma_{\zeta_{2}}-1)\cdots(\gamma_{\zeta_{i-\nu-1}}-1)]_{l},\\
\zeta_{1}<\zeta_{2}<\dots<\zeta_{i-\nu-1},\zeta_{s}\in\{1,2,\dots,i-1\}, 
\end{split}
\end{equation}
where $[\cdot]_{l}$ denotes the $l^{th}$ combination of the product of the elements in parenthesis, therefore we have $l_{\text{max}}=\binom{i-1}{i-\nu-1}$. $\zeta_{s}$ denote all $\zeta$ in \eqref{eq:z-value}.
For the case $\nu=i-1$, if $2\le i$, we define $Z_{\nu,i}=-1$; if $i=1$, we define $Z_{\nu,i}=1$. Beside that, we define $Z_{\nu,i}=0$ for the case $\nu=i$. 
\begin{remark}
\label{rem: different-relax-techniques}
The decay rate in \eqref{eq:convex-hocbf-constraint} used by the iMPC-DHOCBF is partially relaxed compared to the one in \eqref{eq:nonconvex-hocbf-constraint} due to the requirement of the linearization. This can affect the feasibility of the optimization.
\end{remark}

\subsection{CFTOC Problem}
\label{subsec:convex-mpc-dhocbf}

In Secs.~\ref{subsec:linearization-dynamics} and ~\ref{subsec:linearization-dhocbf}, we have illustrated the linearization of system dynamics as well as the safety constraints with DHOCBF.
This allows us to consider them as constraints into a convex MPC formulation at each iteration, which we call convex finite-time constrained optimization control (CFTOC). This is solved at iteration $j$ with optimization variables $\mathbf{U}_{t}^{j} = [ \mathbf{u}_{t,0}^{j}, \dots, \mathbf{u}_{t,N-1}^{j}]$ and $\Omega_{t, i}^{j} = [\omega_{t,0,i}^{j},\dots, \omega_{t,N,i}^{j}]$, where $i\in \{1,\dots, m_{\text{cbf}}\}$.

\noindent\rule{\columnwidth}{0.4pt}
  \textbf{CFTOC of iMPC-DHOCBF at iteration $j$:}  
\begin{subequations}
{\small
\label{eq:impc-dcbf}
\begin{align}
\label{eq:impc-dcbf-cost}
  & \min_{\mathbf{U}_{t}^{j},\Omega_{t,1}^{j},\dots, \Omega_{t,m_{\text{cbf}}}^{j}} p(\mathbf{x}_{t,N}^{j})+\sum_{k=0}^{N-1} q(\mathbf{x}_{t,k}^{j},\mathbf{u}_{t,k}^{j},\omega_{t,k,i}^{j}) \\
   \text{s.t.} \ & \mathbf{x}_{t,k+1}^{j}{-} \bar{\mathbf{x}}_{t,k+1}^{j}{=}A^{j}(\mathbf{x}_{t,k}^{j}-\bar{\mathbf{x}}_{t,k}^{j}){+}B^{j}(\mathbf{u}_{t,k}^{j}-\bar{\mathbf{u}}_{t,k}^{j}), \label{subeq:impc-dcbf-linearized-dynamics}\\
    & \mathbf{u}_{t,k}^{j} \in \mathcal U, \  \mathbf{x}_{t,k}^{j} \in \mathcal X, \ \omega_{t,k,i}^{j}\in \mathbb{R},\label{subeq:impc-dcbf-variables-bounds}\\
    & \tilde{\psi}_{i-1}(\mathbf{x}_{t,k}^{j})+ \sum_{\nu=1}^{i}Z_{\nu,i}  (1-\gamma_{i})^{k}\tilde{\psi}_{0}(\mathbf{x}_{t,\nu}^{j})  \ge \nonumber \\ 
    & \omega_{t,k,i}^{j}Z_{0,i}(1-\gamma_{i})^{k}\tilde{\psi}_{0}(\mathbf{x}_{t,0}^{j}) \label{subeq:impc-dcbf-linearized-hocbf},
\end{align}
}
\end{subequations}
\noindent\rule{\columnwidth}{0.4pt}
In the CFTOC, the linearized dynamics constraints in~\eqref{eq:linearized-dynamics} and the linearized DHOCBF constraints in~\eqref{eq:convex-hocbf-constraint} are enforced with constraints~\eqref{subeq:impc-dcbf-linearized-dynamics} and ~\eqref{subeq:impc-dcbf-linearized-hocbf} at each open loop time step $k\in\{0,\dots, N-1\}$.
The state and input constraints are considered in~\eqref{subeq:impc-dcbf-variables-bounds}. 
The slack variables are unconstrained as the goal of the optimization itself is to minimize the deviation from the nominal DHOCBF constraints with cost term $q(\cdot, \cdot, \omega_{t,k,i}^j)$, while ensuring feasibility of the optimization, as discussed in~\cite{zeng2021decay}.
Note that, for ensuring the safety guarantee established by the DHOCBF, the constraints~\eqref{subeq:impc-dcbf-linearized-hocbf} are enforced with $i\in\{0,\dots, m_{\text{cbf}}\}$, where $Z_{\nu,i}\in \mathbb{R}$ is as defined in~\eqref{eq:convex-hocbf-constraint} with $\nu \in\{0,..,i\}$.
The optimal decision variables of \eqref{eq:impc-dcbf} at iteration $j$ is a list of control input vectors as $\mathbf{U}_{t}^{*,j}=[\mathbf{u}_{t,0}^{*,j},\dots,\mathbf{u}_{t,N-1}^{*,j}]$ and a list of slack variable vectors as $\Omega_{t,i}^{*,j}=[\omega_{t,0,i}^{*,j},\dots,\omega_{t,N-1,i}^{*,j}]$.
The CFTOC is solved iteratively in our proposed iMPC-DHOCBF and the solution can be extracted once the convergence criteria or the maximum iteration number $j_{\text{max}}$ is reached, as shown in Alg.~\ref{alg:iMPC-DCBF}.

%% file: sections/results.tex
\begin{figure*}[t]
    \centering
    \begin{subfigure}[t]{0.24\linewidth}
        \centering
        \includegraphics[width=1\linewidth]{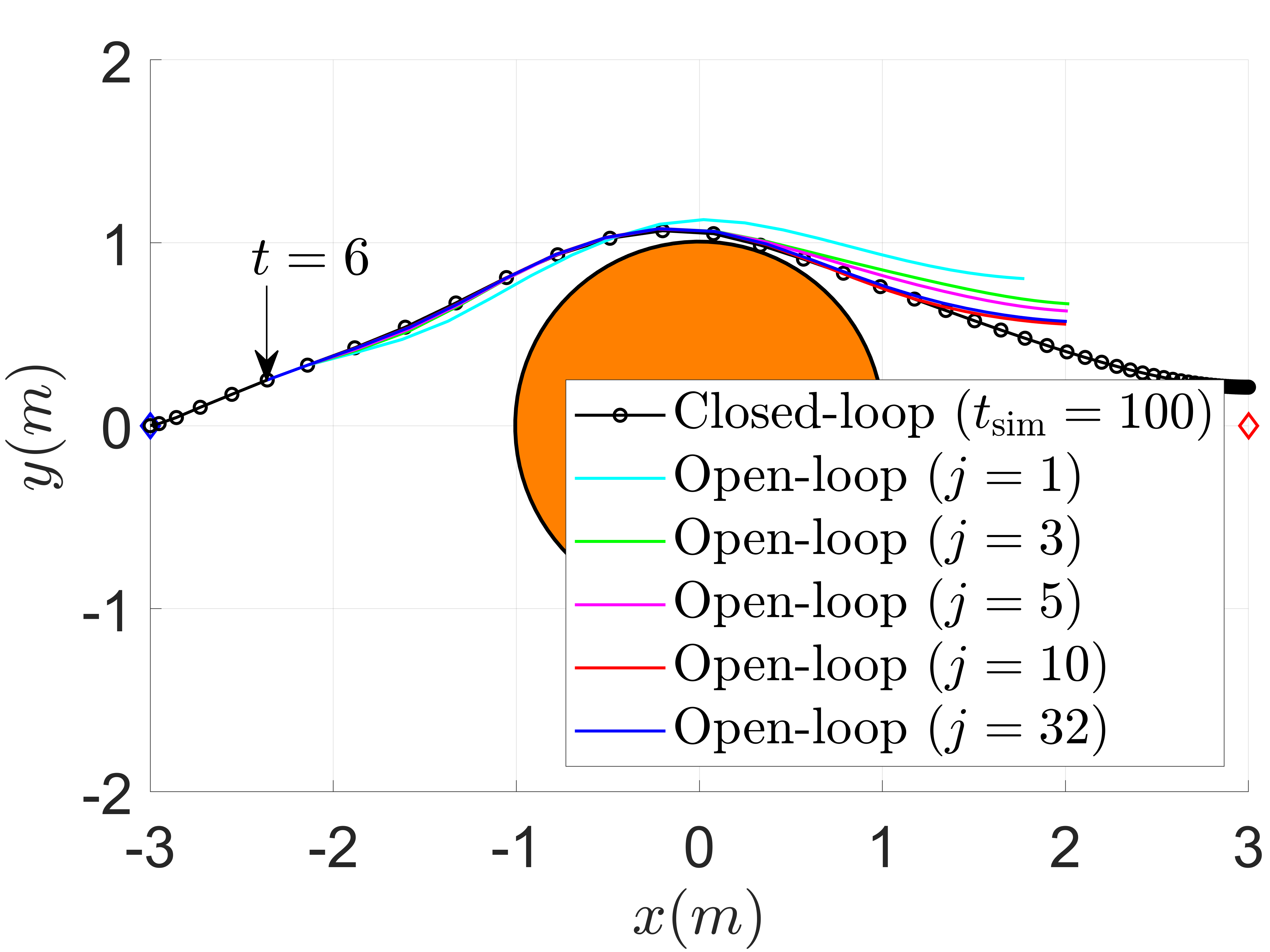}
        \caption{iMPC-DHOCBF when $N=24$, $\gamma_1=\gamma_2 = 0.4$.}
        \label{fig:openloop-snapshots}
    \end{subfigure}
    \begin{subfigure}[t]{0.24\linewidth}
        \centering
        \includegraphics[width=1\linewidth]{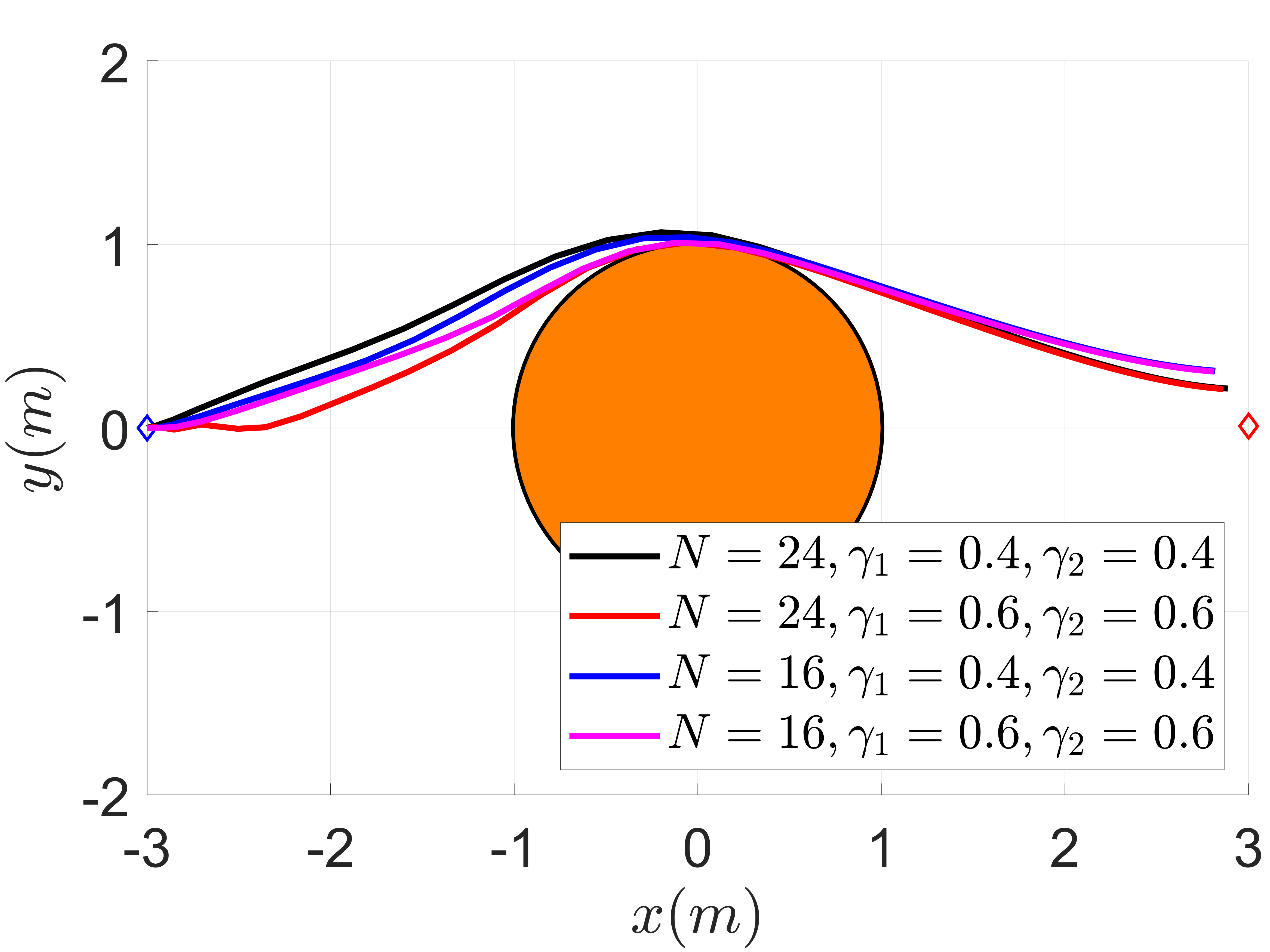}
        \caption{iMPC-DHOCBF with $m_{\text{cbf}} = 2$.}
        \label{fig:closedloop-snapshots1}
    \end{subfigure}  
    \begin{subfigure}[t]{0.24\linewidth}
        \centering
        \includegraphics[width=1\linewidth]{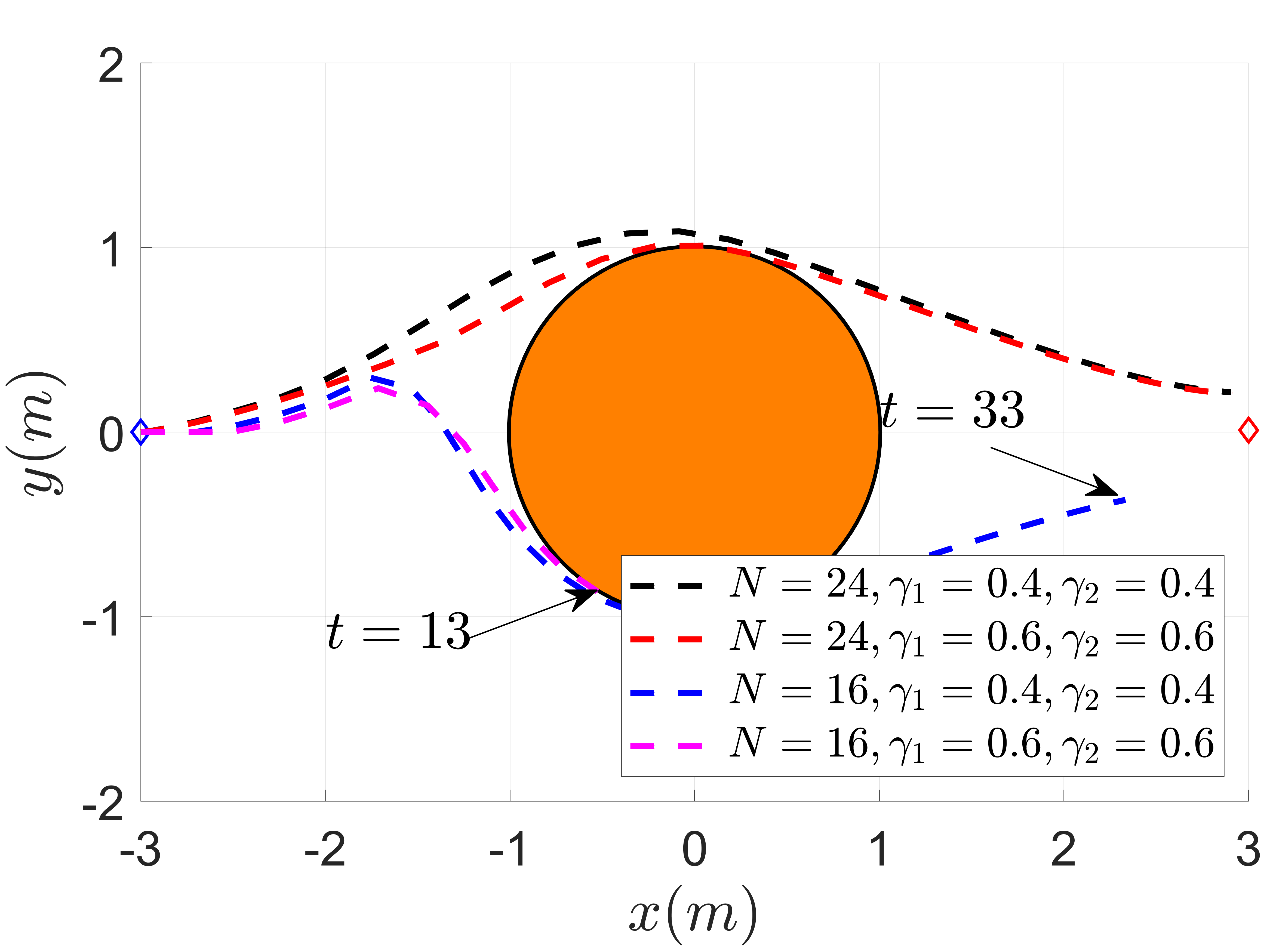}
        \caption{NMPC-DHOCBF with $m_{\text{cbf}} = 2$.}
        \label{fig:closedloop-snapshots2}
    \end{subfigure}
    \begin{subfigure}[t]{0.24\linewidth}
        \centering
        \includegraphics[width=1\linewidth]{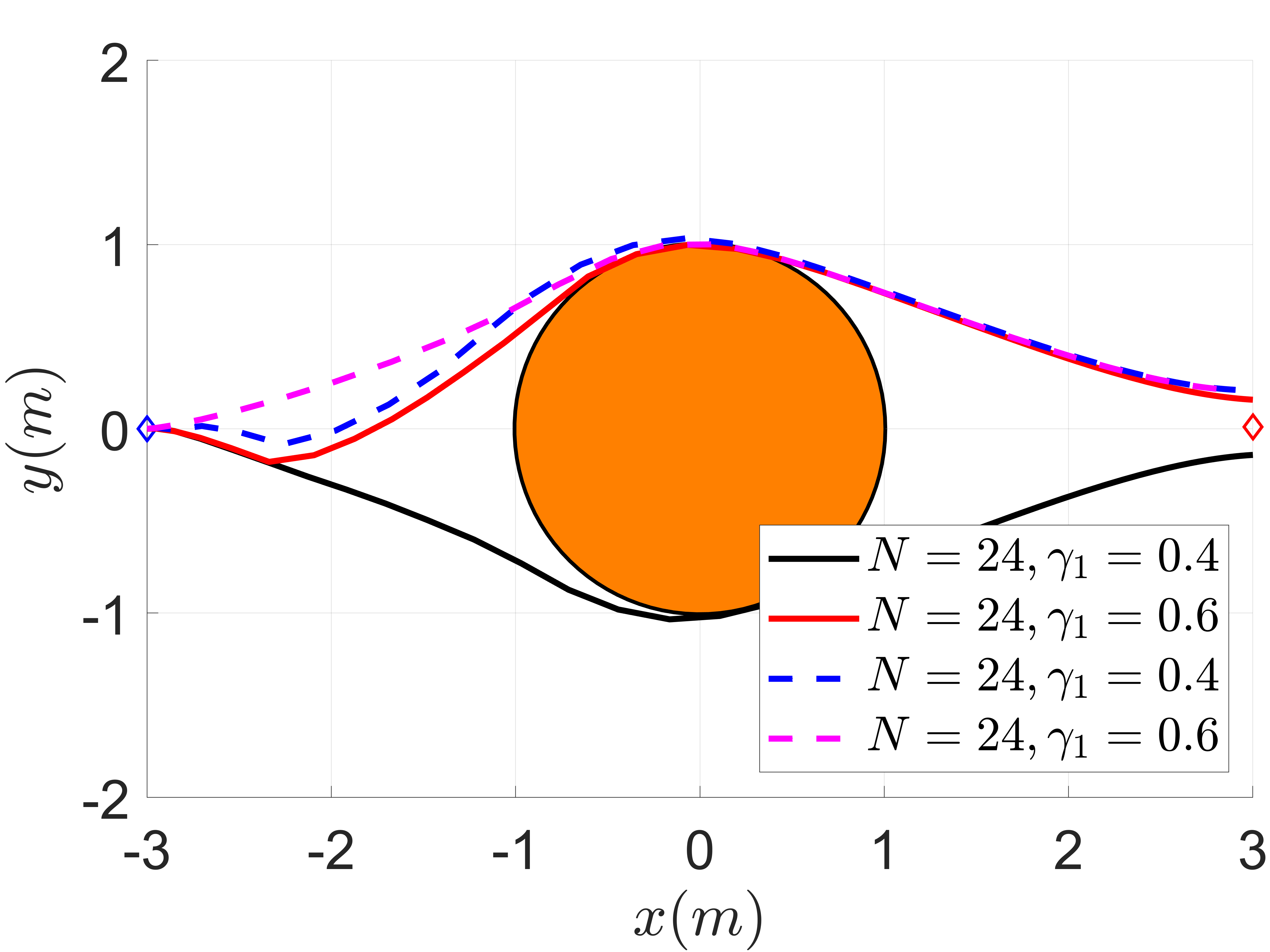}
        \caption{iMPC-DHOCBF and NMPC-DHOCBF with $m_{\text{cbf}}=1$.}
        \label{fig:closedloop-snapshots3}
    \end{subfigure}
    \caption{Open-loop and closed-loop trajectories with controllers iMPC-DHOCBF (solid lines) and NMPC-DHOCBF (dashed lines): (a) several open-loop trajectories at different iterations predicted at $t=6$ and one closed-loop trajectory with controller iMPC-DHOCBF; (b) closed-loop trajectories with controller iMPC-DHOCBF with different choices of $N$ and $\gamma$; (c) closed-loop trajectories with controller NMPC-DHOCBF with different choices of $N$ and $\gamma$. Note that two trajectories stop at $t=13$ and $t=33$ because of infeasibility; (d) closed-loop trajectories with controllers iMPC-DHOCBF and NMPC-DHOCBF with $m_{\text{cbf}}=1$. Both methods work well for safety-critical navigation.
    } 
    \label{fig:open-closed-loop}
\end{figure*}

\begin{figure*}[t]
    \centering
    \begin{subfigure}[t]{0.24\linewidth}
        \centering
        \includegraphics[width=1.0\linewidth]{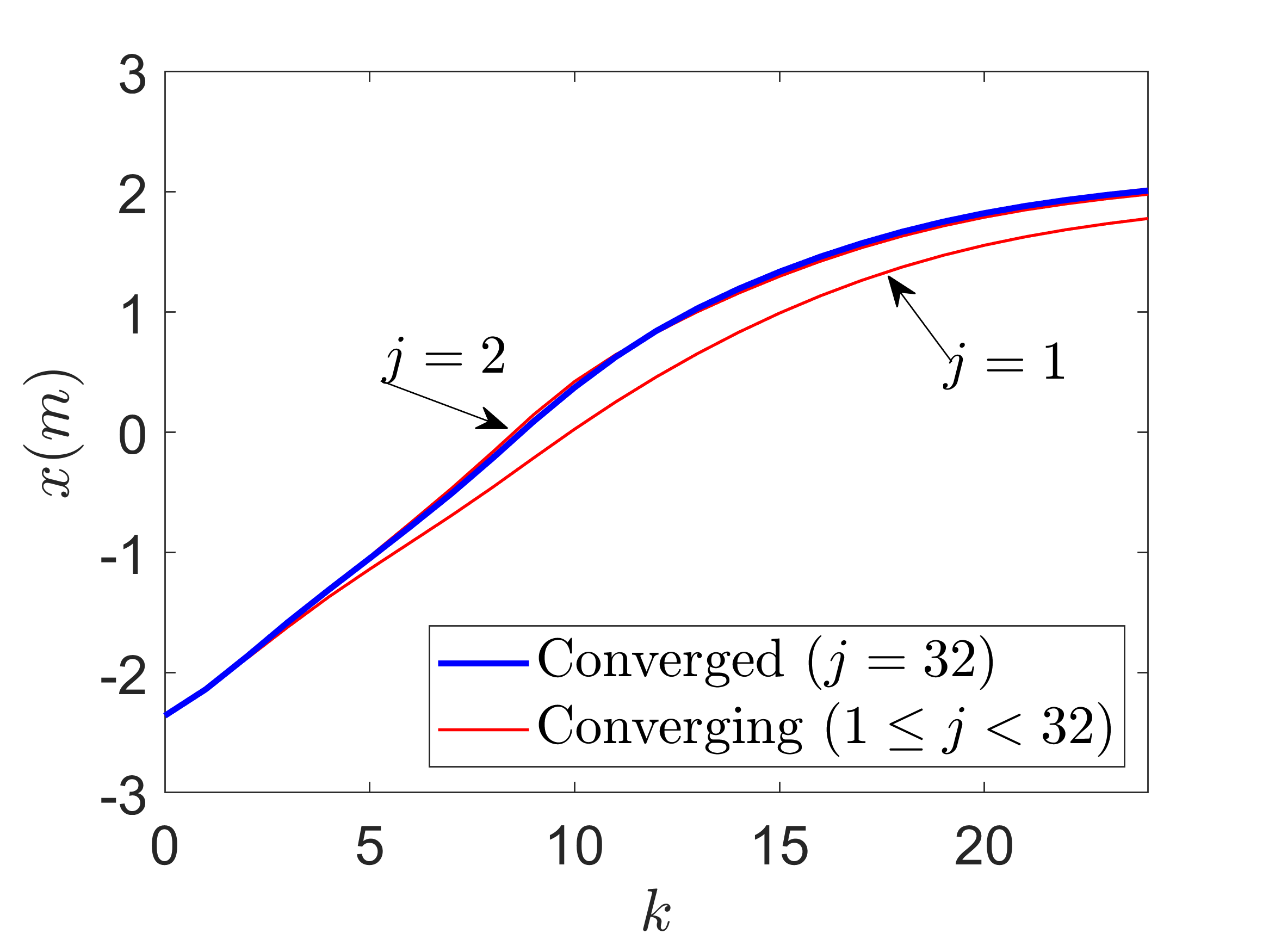}
        \caption{Location $x$}
        \label{subfig:convergence-state-1}
    \end{subfigure}
    \begin{subfigure}[t]{0.24\linewidth}
        \centering
        \includegraphics[width=1.0\linewidth]{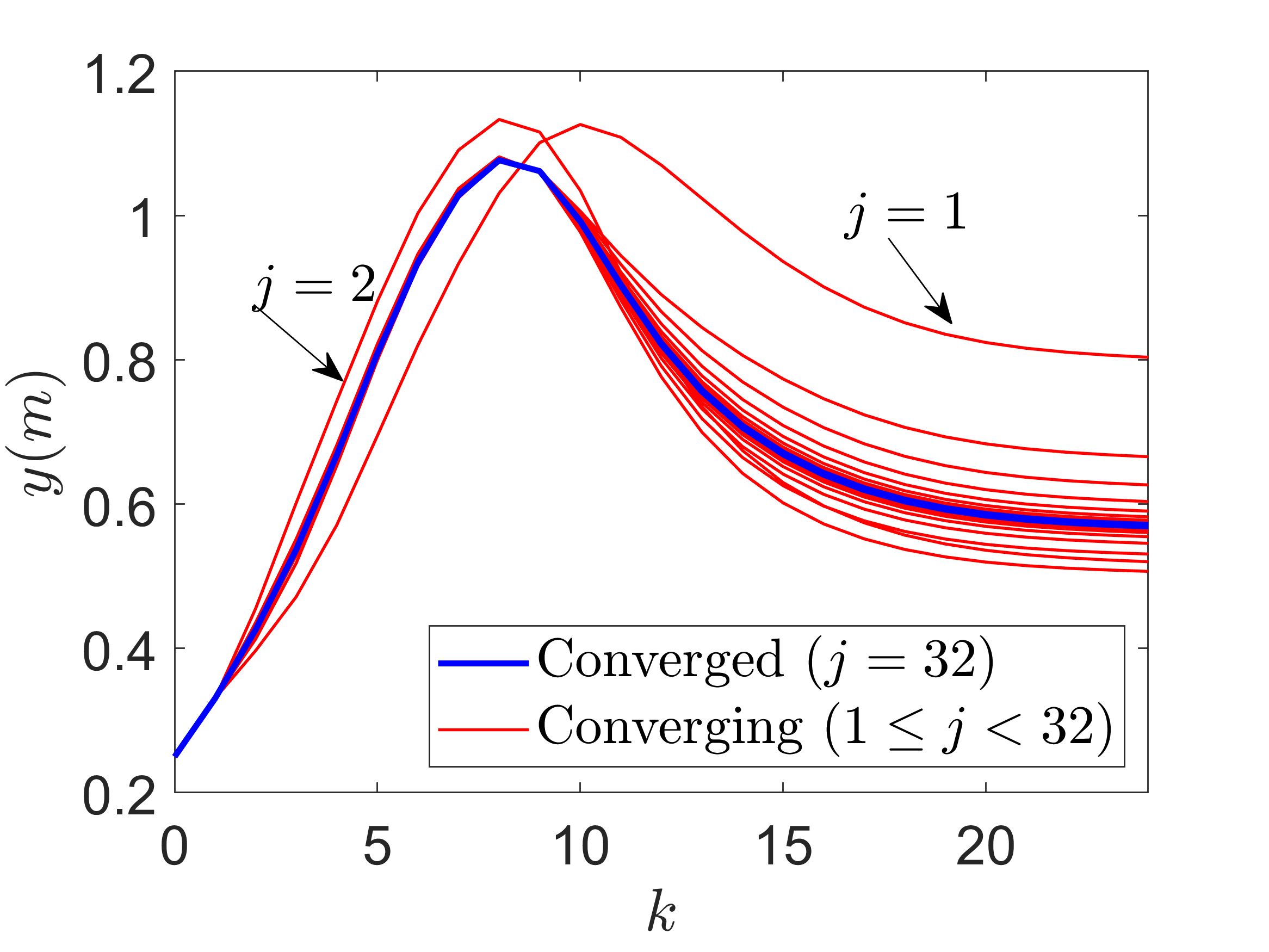}
        \caption{Location $y$}
        \label{subfig:convergence-state-2}
    \end{subfigure}  
    \begin{subfigure}[t]{0.24\linewidth}
        \centering
        \includegraphics[width=1.0\linewidth]{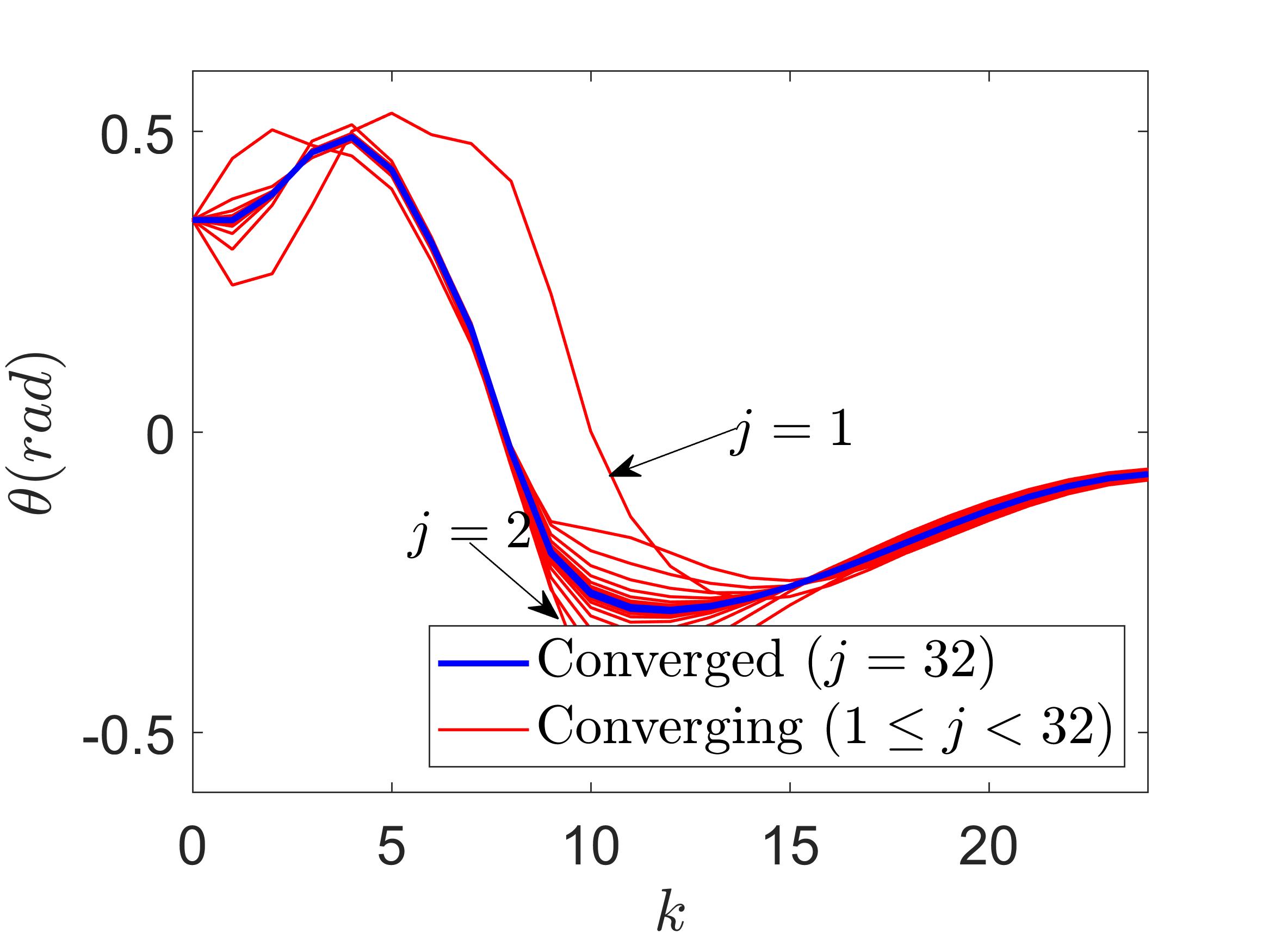}
        \caption{Orientation $\theta$}
        \label{subfig:convergence-state-3}
    \end{subfigure}
    \begin{subfigure}[t]{0.24\linewidth}
        \centering
        \includegraphics[width=1.0\linewidth]{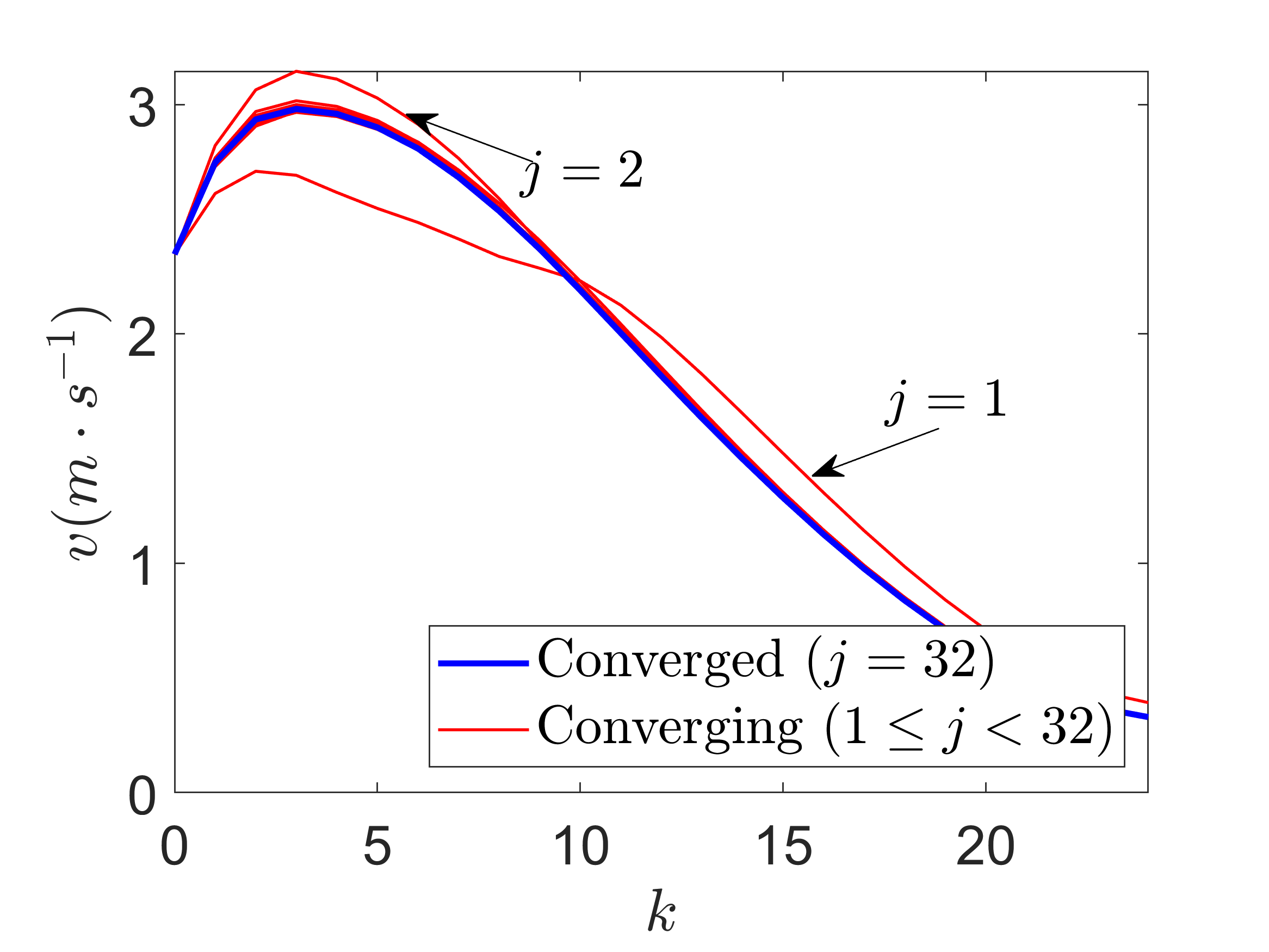}
        \caption{Speed $v$}
        \label{subfig:convergence-state-4}
    \end{subfigure}
    \caption{Iterative convergence of all states at converged iteration $j_{6,\text{conv}}=32$ with $N=24, m_{\text{cbf}}=2$, $\gamma_{1}=\gamma_{2}=0.4$. iMPC does help to optimize the cost function to reach local optimal minimum.}
    \label{fig:convergence-state}
\end{figure*}

\section{Case Study}

In this section, we present numerical results to validate our proposed approach using a unicycle model. We provide a performance comparison with the baseline NMPC-DHOCBF approach. The NMPC-DHOCBF is simply extended by using relaxed DHOCBF based on~\eqref{subeq:mpc-cbf-cons2} in NMPC-DCBF~\eqref{eq:mpc-cbf}, as discussed in~\cite[Rem. 4]{zeng2021enhancing}.

\subsection{Numerical Setup}

\subsubsection{System Dynamics}
Consider a discrete-time unicycle model in the form
\begin{equation}
\label{eq:unicycle-model}
\begin{bmatrix} x_{t+1}{-}x_t \\ y_{t+1}{-}y_t \\ \theta_{t+1}{-}\theta_t \\ v_{t+1}{-}v_t \end{bmatrix}{=}\begin{bmatrix} v_{t} \cos(\theta_{t}) \Delta t \\ v_{t}\sin(\theta_{t}) \Delta t \\ 0 \\ 0 \end{bmatrix}{+}\begin{bmatrix} 0 & 0 \\ 0 & 0 \\ \Delta t & 0 \\ 0 & \Delta t \end{bmatrix}
\begin{bmatrix} u_{1,t} \\ u_{2,t} \end{bmatrix},
\end{equation}
where $\mathbf{x}_{t}=[x_{t},y_{t},\theta_{t},v_{t}]^{T}$ captures the 2-D location, heading angle, and linear speed; $\mathbf{u}_{t}=[u_{1,t},u_{2,t}]^{T}$ represents angular velocity ($u_{1}$) and linear acceleration ($u_{2}$), respectively.
The system is discretized with $\Delta t = 0.1$.
System~\eqref{eq:unicycle-model} is subject to the following state and input constraints:
\begin{equation}
\begin{split}
\label{eq:state-input-constraint}
\mathcal{X}&=\{\mathbf{x}_{t}\in \mathbb{R}^{4}: -10\cdot \mathcal{I}_{4\times1} \le \mathbf{x}_{t}\le 10\cdot \mathcal{I}_{4\times1}\},\\
\mathcal{U}&=\{\mathbf{u}_{t}\in \mathbb{R}^{2}: [-7, -5]^{T} \le \mathbf{u}_{t}\le [7, 5]^{T}\}.
\end{split}
\end{equation}

\subsubsection{System Configuration}
The initial state is $[-3,0,0,0]^{T}$ and the target state is $\mathbf{x}_{r}=[3,0.01,0,0]^{T},$ which are marked as blue and red diamonds in Fig.~\ref{fig:open-closed-loop}. The circular obstacle is centered at $(0,0)$ with $r=1,$, which is displayed in orange. The other reference vectors are $\mathbf{u}_{r}=[0,0]^{T}$ and $\omega_{r}=[1,1]^{T}$. 
We use the offset $y=0.01m$ in $\mathbf{x}_{r}$ to prevent singularity of the optimization problem. 

\subsubsection{DHOCBF}
As a candidate DHOCBF function $\psi_{0}(\mathbf{x}_{t})$, we choose a quadratic distance function for circular obstacle avoidance $h(\mathbf{x}_{t})= (x_{t}-x_{0})^{2}+(y_{t}-y_{0})^{2}-r^{2}$, where $(x_{0},y_{0})$ and $r$ denote the obstacle center location and radius, respectively. The linearized DHOCBF $\tilde{\psi}_{0}(\mathbf{x}_{t,k}^{j})$ in \eqref{eq:linearized-CBFs} is defined as $\tilde{\psi}_{0}(\mathbf{x}_{t,k}^{j}) \coloneqq h_{\parallel}(\mathbf{x}_{t,k}^{j}|\tilde{\mathbf{x}}_{t,k}^{j})$, with
\begin{equation}
\begin{split}
h_{\parallel}(\mathbf{x}_{t,k}^{j}|\tilde{\mathbf{x}}_{t,k}^{j})=(\tilde{x}_{t,k}^{j}-x_{0})x_{t,k}^{j}+(\tilde{y}_{t,k}^{j}-y_{0})y_{t,k}^{j}\\
-(r^{2}-x_{0}^{2}-y_{0}^2+\tilde{x}_{t,k}^{j}x_{0}+\tilde{y}_{t,k}^{j}y_{0}), 
\end{split}
\end{equation}
where $h_{\parallel}(\mathbf{x}_{t,k}^{j}|\tilde{\mathbf{x}}_{t,k}^{j})$ is the linearized boundary, whose relative degree is 2; $(\tilde{x}_{t,k}^{j}, \tilde{y}_{t,k}^{j})$ denotes the tangent point of the circular boundary $h(\mathbf{x}_{t})$. From \eqref{eq:z-value}, we have $Z_{0,2}=\gamma_{1}-1,\ Z_{1,2}=-1,\ Z_{0,1}=1,\ Z_{2,2}=Z_{1,1}=0$.

\subsubsection{MPC Design}
The cost function of the MPC problem consists of stage cost
$q(\mathbf{x}_{t,k}^j,\mathbf{u}_{t,k}^{j},\omega_{t,k}^{j})= \sum_{k=0}^{N-1} (||\mathbf{x}_{t,k}^{j}-\mathbf{x}_{r}||_Q^2 + ||\mathbf{u}_{t,k}^{j}-\mathbf{u}_{r}||_R^2 +||\omega_{t,k}^{j}-\omega_{r}||_S^2)$
and terminal cost $p(\mathbf{x}_{t,N}^{j})=||\mathbf{x}_{t,N}^{j}-\mathbf{x}_{r}||_P^2$, where $Q=P=10\cdot \mathcal{I}_{4}, R= \mathcal{I}_{2}$ and $S=1000\cdot \mathcal{I}_{2}$.

\subsubsection{Convergence Criteria}
We use the following absolute and relative convergence functions as convergence criteria mentioned in Alg. \ref{alg:iMPC-DCBF}:
\begin{equation}
\label{eq:convergence-criteria}
    \begin{split}
    e_{\text{abs}}(\mathbf{X}_{t}^{*, j}, \mathbf{U}_{t}^{*, j}) &= ||\mathbf{X}^{*, j} - \bar{\mathbf{X}}^{*, j}|| \\
    e_{\text{rel}}(\mathbf{X}_{t}^{*, j}, \mathbf{U}_{t}^{*, j}, \bar{\mathbf{X}}_{t}^{j}, \bar{\mathbf{U}}_{t}^{j}) &= ||\mathbf{X}^{*, j} - \bar{\mathbf{X}}^{*, j}||/||\bar{\mathbf{X}}^{*, j}||.
\end{split}
\end{equation}
The iterative optimization stops when $e_{\text{abs}} < \varepsilon_{\text{abs}}$ or $e_{\text{rel}} < \varepsilon_{\text{rel}}$, where $\varepsilon_{\text{abs}} = 10^{-4}$, $\varepsilon_{\text{rel}} = 10^{-2}$ and the maximum iteration number is set as $j_{\text{max}} = 1000$.

To make a fair comparison with NMPC-DHOCBF, the hyperparameters $P, Q, R, S$ remain unchanged for all setups.

\subsubsection{Solver Configurations and CPU Specs.}
For iMPC-DHOCBF, we used OSQP~\cite{stellato2020osqp} to solve the convex optimizations at all iterations.
The baseline approach NMPC-DHOCBF is open-source, and was solved using IPOPT~\cite{biegler2009large} with the modeling language Yalmip~\cite{lofberg2004yalmip}.
We used a Windows desktop with Intel Core i7-8700 (CPU 3.2 GHz) running Matlab for all computations.

\subsection{Performance}

\subsubsection{Iterative Convergence}
\label{subsubsec:iterative-convergence}
\begin{figure*}[t]
    \centering
    \begin{subfigure}[t]{0.24\linewidth}
        \centering
        \includegraphics[width=1.0\linewidth]{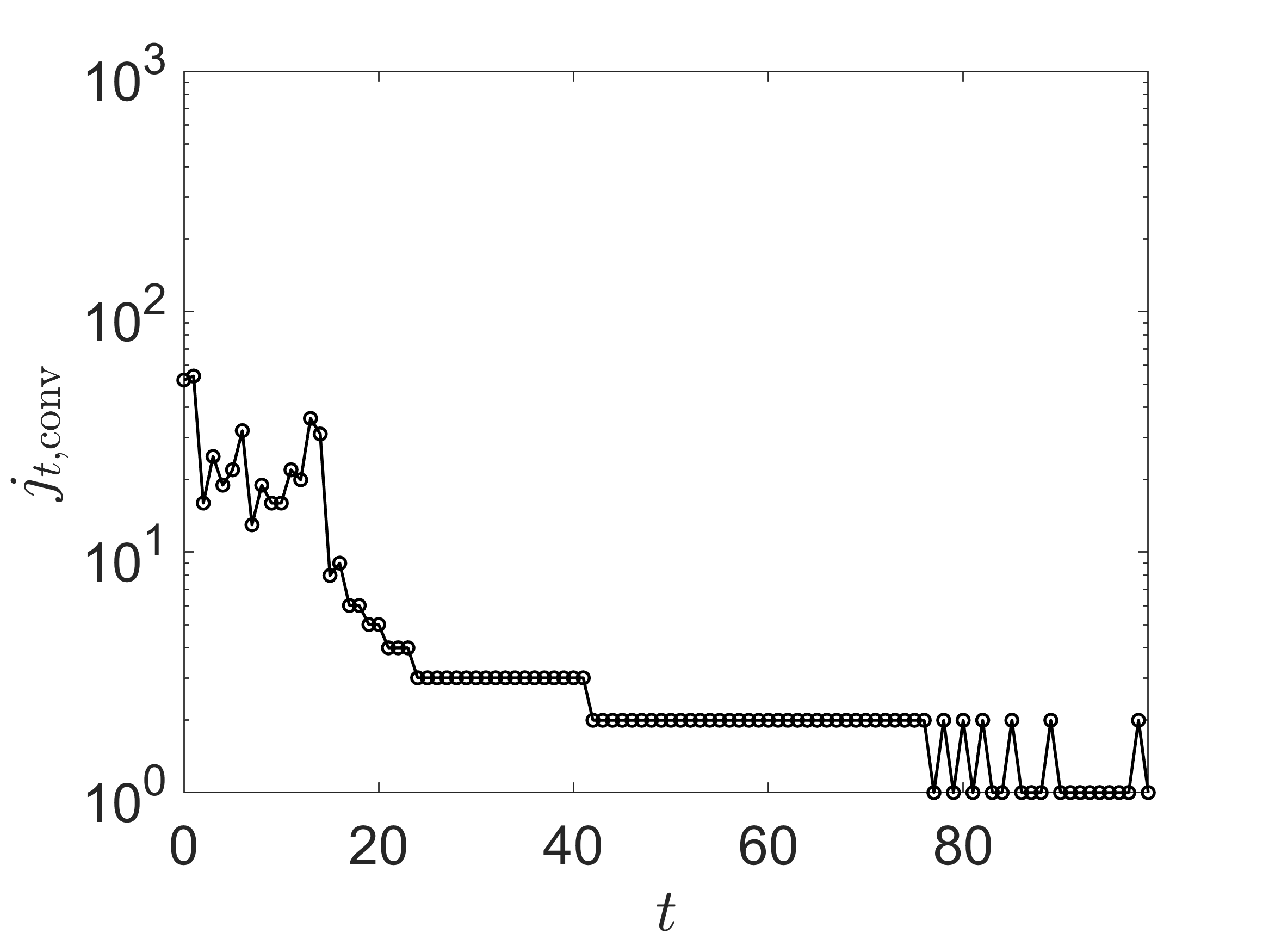}
        \caption{$\gamma_1 = 0.4, \gamma_2 = 0.4$}
        \label{subfig:convergence-J-function-1}
    \end{subfigure}
    \begin{subfigure}[t]{0.24\linewidth}
        \centering
        \includegraphics[width=1.0\linewidth]{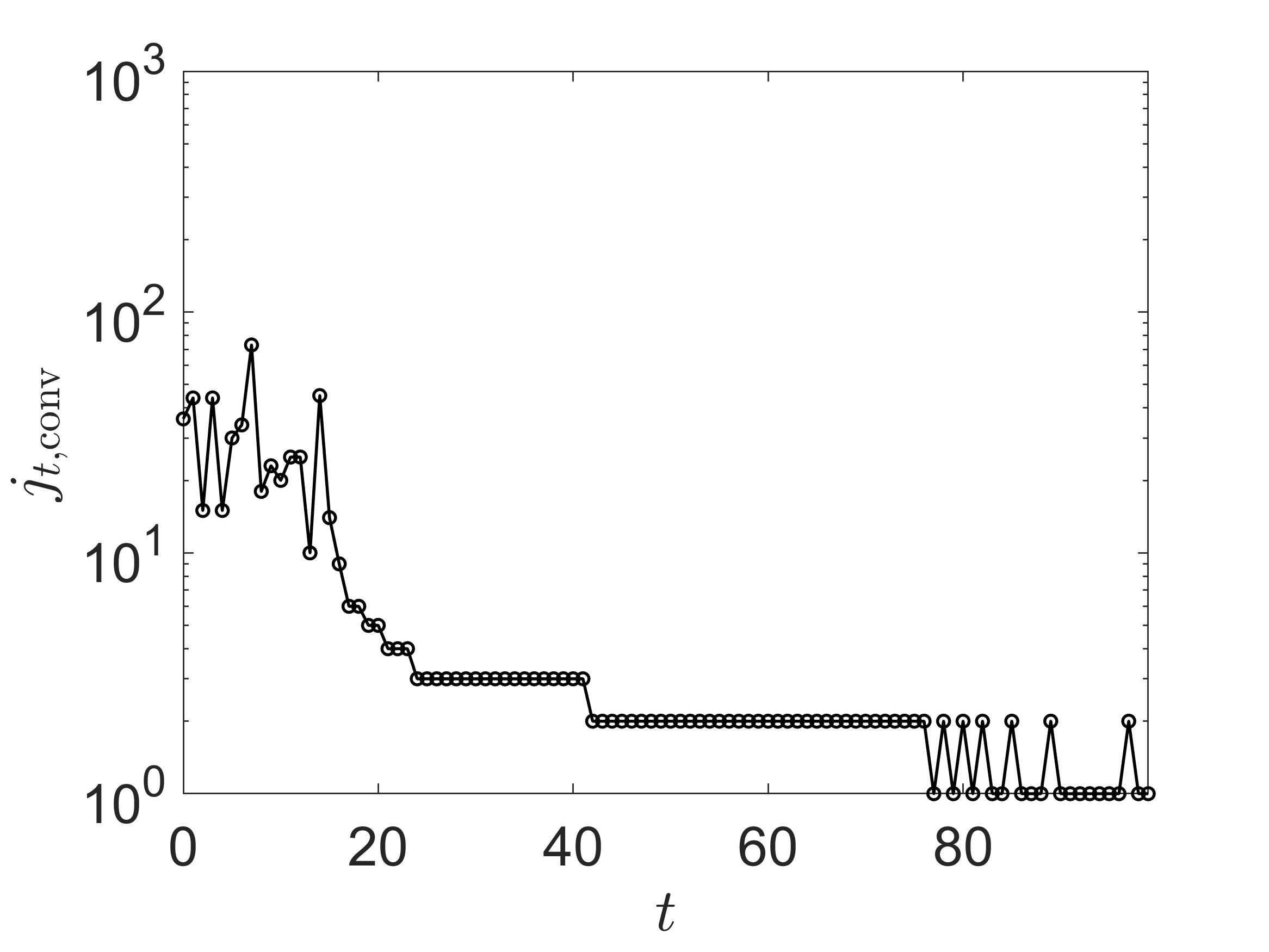}
        \caption{$\gamma_1 = 0.4, \gamma_2 = 0.6$}
        \label{subfig:convergence-J-function-2}
    \end{subfigure}  
    \begin{subfigure}[t]{0.24\linewidth}
        \centering
        \includegraphics[width=1.0\linewidth]{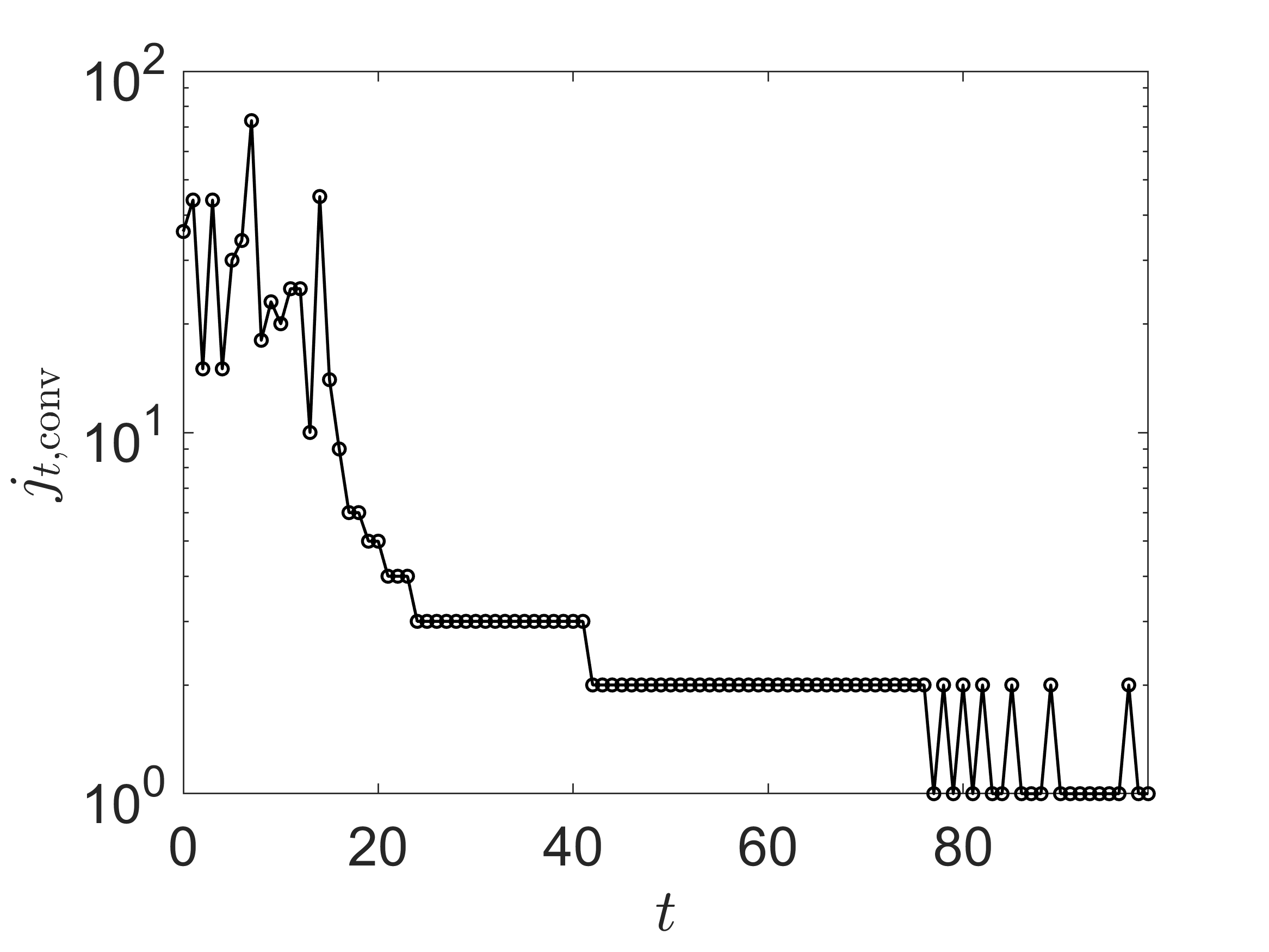}
        \caption{$\gamma_1 = 0.6, \gamma_2 = 0.4$}
        \label{subfig:convergence-J-function-3}
    \end{subfigure}
    \begin{subfigure}[t]{0.24\linewidth}
        \centering
        \includegraphics[width=1.0\linewidth]{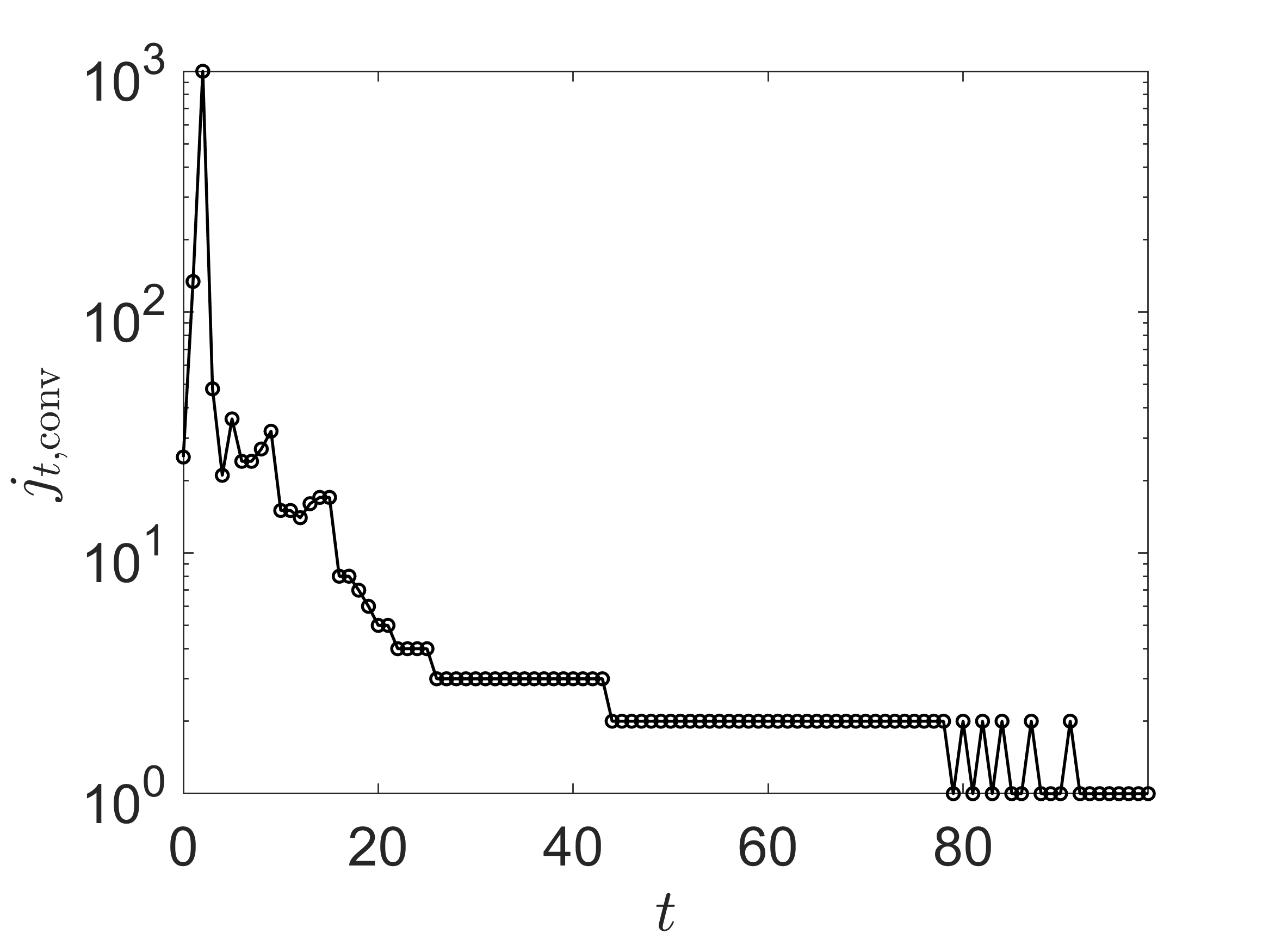}
        \caption{$\gamma_1 = 0.6, \gamma_2 = 0.6$}
        \label{subfig:convergence-J-function-4}
    \end{subfigure}
    \caption{Number of iterations $j_{t,\text{conv}}$ at each time-step using controller iMPC-DHOCBF with different values of hyperparameters $\gamma_1, \gamma_2$ with $N=24, m_{\text{cbf}}=2$.
    We can observe that, for almost all time-steps, the iterative optimization converges within 100 iteraitons ($j_{t, \text{conv}} < 10^2$), which is affected very little with respect to hyperparameters.}
    \label{fig:convergence-J-function}
\end{figure*}

The iterative convergence is shown in Figs. \ref{fig:openloop-snapshots}, \ref{fig:convergence-J-function} and \ref{fig:convergence-state}.
Fig. \ref{fig:openloop-snapshots} shows the closed-loop trajectory (the black line) generated by solving the iMPC-DHOCBF until the converged iteration $j_{t,\text{conv}}$ from $t=0$ to $t=t_{\text{sim}}=100$ and open-loop iteratively converging trajectories (colored lines) at different iterations at $t=6$.
Fig.~\ref{fig:convergence-state} presents more details on the iterative convergence of states at different iterations at $t=6$ with number of iterations $j_{t,\text{conv}} = 32$.
We note that, after around 10 iterations, the converging lines for the states (red lines) nearly overlap with the converged line (blue line) in Fig. \ref{fig:convergence-state}. This verifies the relations of the converging trajectory (red line) and the converged trajectory (blue line) in Fig. \ref{fig:openloop-snapshots}.
The optimization is shown to converge at iteration $j_{t,\text{conv}}$ at time step $t$ for different hyperparameters $\gamma$ under specific convergence criteria~\eqref{eq:convergence-criteria}, shown in Fig. \ref{fig:convergence-J-function}.
We can see that for the first 15 time steps the iMPC-DHOCBF triggers more iterations to drive the system to avoid the obstacle than time steps after 20 where the system already passes the obstacle.
The maximum converged iteration $j_{t,\text{conv}}$ is 1000 at time step $t=2$ in Fig. \ref{subfig:convergence-J-function-4} with $\gamma_{1}=\gamma_{2}=0.6$, which reveals that the peak of the converged iteration over time increases if we choose larger $\gamma$.
For the majority of the time-steps, the iterative optimization converges within 100 iterations ($j_{t, \text{conv}} < 100$).

\begin{table*}[]
\centering
\resizebox{0.99\textwidth}{!}{
\begin{tabular}{|cc|cccccc|}
\hline
\multicolumn{2}{|c|}{Approaches} & $N = 4$ & $N = 8$ & $N = 12$ & $N = 16$ & $N = 20$ & $N = 24$ \\ \hline
\multicolumn{1}{|c|}{\multirow{2}{*}{\begin{tabular}[c]{@{}c@{}}NMPC-DHOCBF\\ ($m_{\text{cbf}} = 2$)\end{tabular}}} & mean / std (s) & $3.687\pm6.360$ & $23.882\pm17.988$ & $27.329\pm20.115$ & $28.953\pm22.058$ & $30.970\pm23.564$ & $29.929\pm22.105$ \\
\multicolumn{1}{|c|}{} & infeas. rate & 5.8\% & 27.5\% & 21.1\%  &  16.4\% & 14.5\% & 14.4\% \\ \hline
\multicolumn{1}{|c|}{\multirow{2}{*}{\begin{tabular}[c]{@{}c@{}}NMPC-DHOCBF\\ ($m_{\text{cbf}}  = 1$)\end{tabular}}} & mean / std (s) & $2.933\pm4.678$ & $19.077\pm14.024$ & $20.418\pm15.401$ & $22.749\pm17.039$ & $24.053\pm17.811$ & $25.365\pm18.211$ \\
\multicolumn{1}{|c|}{} & infeas. rate & 6.3\% & 13.9\% & 13.0\% & 14.6\% & 13.8\% & 15.4\% \\ \hline
\multicolumn{1}{|c|}{\multirow{2}{*}{\begin{tabular}[c]{@{}c@{}}iMPC-DHOCBF\\ ($m_{\text{cbf}}  = 2$)\end{tabular}}} & mean / std (s) & $0.135\pm0.294$ & $0.104\pm0.242$ & $0.102\pm0.217$ & $0.131\pm0.301$ & $0.165\pm0.400$ & $0.135\pm0.274$ \\
\multicolumn{1}{|c|}{} & infeas. rate & 6.3\% & 8.0\% & 10.4\% & 10.9\% & 10.9\% & 10.2\% \\ \hline
\multicolumn{1}{|c|}{\multirow{2}{*}{\begin{tabular}[c]{@{}c@{}}iMPC-DHOCBF\\ ($m_{\text{cbf}}  = 1$)\end{tabular}}} & mean / std (s) & $0.131\pm0.286$ & $0.114\pm0.260$ & $0.109\pm0.237$ & $0.137\pm0.316$ & $0.173\pm0.414$ & $0.152\pm0.317$ \\
\multicolumn{1}{|c|}{} & infeas. rate & 6.3\% & 8.0\% & 10.4\% & 10.9\% & 10.9\% & 11.1\% \\ \hline
\end{tabular}
}
\caption{Statistical benchmark for computation time and feasibility between NMPC-DHOCBF and iMPC-DHOCBF with randomized states.
The target position is shared among four approaches and the hyperparameters are fixed as $\gamma_{1}=\gamma_{2}=0.4$ for all random scenarios. 
\label{tab:compuation-time1}}
\end{table*}

\begin{table*}[]
\centering
\resizebox{0.99\textwidth}{!}{
\begin{tabular}{|cc|cccccc|}
\hline
\multicolumn{2}{|c|}{Approaches} & $N = 4$ & $N = 8$ & $N = 12$ & $N = 16$ & $N = 20$ & $N = 24$ \\ \hline
\multicolumn{1}{|c|}{\multirow{2}{*}{\begin{tabular}[c]{@{}c@{}}NMPC-DHOCBF\\ ($m_{\text{cbf}} = 2$)\end{tabular}}} & mean / std (s) & $3.744\pm6.445$ & $28.779\pm20.755$ & $31.319\pm21.921$ & $33.678\pm25.328$ & $36.430\pm26.959$ & $39.543\pm29.941$ \\
\multicolumn{1}{|c|}{} & infeas. rate & 5.6\% & 28.0\% & 20.9\%  &  16.8\% & 17.0\% & 14.6\% \\ \hline
\multicolumn{1}{|c|}{\multirow{2}{*}{\begin{tabular}[c]{@{}c@{}}NMPC-DHOCBF\\ ($m_{\text{cbf}}  = 1$)\end{tabular}}} & mean / std (s) & $3.032\pm4.536$ & $21.414\pm16.518$ & $23.121\pm17.544$ & $24.011\pm17.711$ & $26.599\pm19.480$ & $29.671\pm20.026$ \\
\multicolumn{1}{|c|}{} & infeas. rate & 6.4\% & 17.0\% & 15.2\% & 15.5\% & 16.7\% & 13.2\% \\ \hline
\multicolumn{1}{|c|}{\multirow{2}{*}{\begin{tabular}[c]{@{}c@{}}iMPC-DHOCBF\\ ($m_{\text{cbf}}  = 2$)\end{tabular}}} & mean / std (s) & $0.158\pm0.326$ & $0.134\pm0.279$ & $0.163\pm0.353$ & $0.163\pm0.373$ & $0.184\pm0.398$ & $0.164\pm0.344$ \\
\multicolumn{1}{|c|}{} & infeas. rate & 6.1\% & 8.0\% & 10.2\% & 10.7\% & 10.8\% & 10.8\% \\ \hline
\multicolumn{1}{|c|}{\multirow{2}{*}{\begin{tabular}[c]{@{}c@{}}iMPC-DHOCBF\\ ($m_{\text{cbf}}  = 1$)\end{tabular}}} & mean / std (s) & $0.167\pm0.340$ & $0.139\pm0.291$ & $0.170\pm0.362$ & $0.170\pm0.379$ & $0.201\pm0.435$ & $0.176\pm0.378$ \\
\multicolumn{1}{|c|}{} & infeas. rate & 6.1\% & 8.0\% & 10.2\% & 10.7\% & 10.8\% & 10.8\% \\ \hline
\end{tabular}
}
\caption{Statistical benchmark between NMPC-DHOCBF and iMPC-DHOCBF with the same randomized states as in Tab. \ref{tab:compuation-time1}. 
The target position is shared among four approaches and the hyperparameters are fixed as $\gamma_{1}=\gamma_{2}=0.6$ for all scenarios. Based on Tab. \ref{tab:compuation-time1} and Tab. \ref{tab:compuation-time2} we conclude that iMPC-DHOCBF outperforms NMPC-DHOCBF in computing time and infeasibility rate.
\label{tab:compuation-time2}}
\end{table*}

\subsubsection{Convergence with Different Hyperparameters}
Fig. \ref{fig:closedloop-snapshots1}, \ref{fig:closedloop-snapshots2} and \ref{fig:closedloop-snapshots3} show the closed-loop trajectories generated by solving iMPC-DHOCBF (solid lines) and NMPC-DHOCBF (dashed lines) at converged iteration $j_{t,\text{conv}}$ from $t=0$ to $t=t_{\text{sim}}=45$ with different hyperparameters.
Both controllers show good performance on obstacle avoidance.
Based on black, red, blue and magenta lines with the highest order of CBF constraint $m_{\text{cbf}}=2$ in Fig. \ref{fig:closedloop-snapshots1} and \ref{fig:closedloop-snapshots2}, as $\gamma_{1}, \gamma_{2}$ become smaller, the system tends to turn further away from the obstacle when it is getting closer to obstacle, which indicates a safer control strategy.
From the lines in Fig. \ref{fig:closedloop-snapshots3} where $m_{\text{cbf}}=1$, we can see that the system can still safely navigate around the obstacle, although it turns away from the obstacle later than when having one more CBF constraint in Fig. \ref{fig:closedloop-snapshots1} and \ref{fig:closedloop-snapshots2}, indicating that having CBF constraints up to the relative degree enhances safety.
The blue and magenta dashed lines in Fig. \ref{fig:closedloop-snapshots2} stop at $t=33$ and $t=13$ with $N=16$ as infeasibility happens, which shows that a large horizon is needed to generate complete closed-loop trajectories for some hyperparameters by NMPC-DHOCBF, while iMPC-DHOCBF shows less reliance on selection of horizon since it can generate complete closed-loop trajectories with both $N=16$ and $N=24$, as shown in Fig. \ref{fig:closedloop-snapshots1}.

\subsubsection{Computation Time}
In order to compare computational times between our proposed iMPC-DHOCBF and the baseline NMPC-DHOCBF, 1000 independent randomized safe states are generated in state constraint $\mathcal{X}$ in \eqref{eq:state-input-constraint}.
To make a fair comparison, both approaches use the same $N$ and $m_{\text{cbf}}$ and the computational time and feasibility are evaluated at those randomized sample states.
The distributions of the computation times and infeasibility rates in Tab. \ref{tab:compuation-time1} and Tab. \ref{tab:compuation-time2} correspond to generating one time-step trajectories. 
The mean and standard deviation of computation times increase if the horizon $N$ or $m_{\text{cbf}}$ become larger for NMPC-DHOCBF in Tab. \ref{tab:compuation-time1} and Tab. \ref{tab:compuation-time2}.
Different from NMPC-DHOCBF, the computing time is not heavily influenced by $N$ and $m_{\text{cbf}}$ for iMPC-DHOCBF. Based on the data from the two tables,  we also notice that larger hyperparameter values for $\gamma$ will slightly reduce the computation speed for both methods, which is discussed in Sec. \ref{subsubsec:iterative-convergence} and can be attributed to the rise of converged iteration $j_{t,\text{conv}}$.  Compared to NMPC-DHOCBF, the computing speed of our proposed method is much faster with the improvement in computation time directly proportional to the horizon, $e.g.$,  $100\sim 300$ times faster than the baseline given the chosen hyperparameters.

\subsubsection{Optimization Feasibility}
The rate of infeasibility increases if the horizon $N$ increases or $m_{\text{cbf}}$ is lower for iMPC-DHOCBF. However, these two hyperparameters are shown not to  affect the infeasibility rate of the NMPC-DHOCBF method proportionally.
As the horizon increases, the infeasibility rate of iMPC-DHOCBF outperforms that of NMPC-DHOCBF.
The main reasons for this come from the difference in the convergence criteria and relaxation techniques for CBF constraints, discussed in Rem. \ref{rem: different-relax-techniques}. The NMPC-DHOCBF under IPOPT should have more strict convergence criteria compared to iMPC-DHOCBF, which obviously limits its feasibility if the number of horizon is large.
Besides, NMPC-DHOCBF is equipped with relaxed nonlinear CBF constraints \eqref{eq:nonconvex-hocbf-constraint}, while iMPC-DHOCBF has relaxed linear CBF constraints \eqref{eq:convex-hocbf-constraint}.
The linearization of the CBF constraints reduces the feasibility region in the state space, as illustrated in Fig.~\ref{fig:linearization-dhocbf}. This meets the expectation of slight decreased feasibility rate of iMPC-DHOCBF when number of horizon is small.
However, we can see that the decline in feasibility rate due to relaxed technique is noticeably outperformed by flexible convergence criteria  with larger number of horizon $N$ in Tab.~\ref{tab:compuation-time1} and ~\ref{tab:compuation-time2}, which validates our linearization technique in the iterative optimization.

%% file: sections/conclusion.tex
\section{Conclusion \& Future Work}
We proposed an iterative convex optimization procedure for safety-critical model predictive control (iMPC) design. Central to our approach are relaxations for the system dynamics and for discrete time high-order control barrier functions (DHOCBF) in the form of 
linearized constraints. We validated the proposed iMPC-DHOCBF approach by applying it to a model of unicycle navigating in an environment with obstacles. We noticed that the computation times for the iMPC-DHOCBF method significantly outperform the ones corresponding to the baseline, usually with even higher feasibility rate.
There are still some limitations of iMPC-DHOCBF that could be ameliorated. One limitation of the proposed method is its linearly relaxed technique will slightly increase infeasibility rate with small size of the horizon.
Another limitation is that the feasibility of the optimization and system safety are not always guaranteed at the same time in the whole state space. We will address these limitations in future work with better linearization, different relaxed techniques as well as adaptive warm-up and convergence criterion for the optimization problem.